\documentclass[12pt,a4paper]{article}

\usepackage{amsmath}
\usepackage{kmath}
\usepackage{kerkis} 

\usepackage[english]{babel}
\usepackage{amstext}
\usepackage{amssymb}
\usepackage{graphicx}
\usepackage{dsfont}
\usepackage{float}
\usepackage{xcolor}
\usepackage{textcomp}
\usepackage{listings}
\usepackage{relsize}
\usepackage{multicol}
\usepackage{subcaption}
\usepackage{physics}
\usepackage{physics}
\usepackage{color}
\usepackage[scale=0.8]{geometry}

\begin{document}

\begin{center}
	
	{{{\Large\bf Periodically Forced Nonlinear Oscillators With Hysteretic Damping}}}\\
	
	\vspace{0.5cm}
	{\large Anastasios Bountis\footnote{E-mail address: anastasios.bountis@nu.edu.kz}$^{,a}$, Konstantinos Kaloudis\footnote{E-mail address: konst.kaloudis@gmail.com}$^{,b}$ and Christos Spitas\footnote{E-mail address: christos.spitas@nu.edu.kz}$^{,b}$}
	
	\vspace{0.2cm}
	
	$^a$Department of Mathematics, $^b$Department of Mechanical and Aerospace Engineering,\\Nazarbayev University, Kabanbay-Batyr, 53,\\ 010000 Nur-Sultan, Republic of Kazakhstan 
\end{center}
\vspace{0.2in}

\begin{abstract}
\textit{	We perform a detailed study of the dynamics of a nonlinear, one-dimensional oscillator driven by a periodic force under hysteretic damping, whose linear version was originally proposed and analyzed by Bishop in \cite{bishop1955treatment}. We first add a small quadratic stiffness term in the constitutive equation and construct the periodic solution of the problem by a systematic perturbation method, neglecting transient terms as $t\rightarrow \infty$. We then repeat the analysis replacing the quadratic by a cubic term, which does not allow the solutions to escape to infinity. In both cases, we examine the dependence of the amplitude of the periodic solution on the different parameters of the model and discuss the differences with the linear model. We point out certain undesirable features of the solutions, which have also been alluded to in the literature for the linear Bishop's model, but persist in the nonlinear case as well. Finally, we discuss an alternative hysteretic damping oscillator model first proposed by Reid \cite{reid1956free}, which appears to be free from these difficulties and exhibits remarkably rich dynamical properties when extended in the nonlinear regime.}
\end{abstract}

\section{Introduction}

As engineering materials are used in increasingly mission-critical (and thus behavior-critical) contexts, where no longer strength but also noise and vibration behavior matter across a spectrum of operating frequencies, their nonlinear characteristics cannot afford to be ignored \cite{OlejAwre2018}. Power law materials are used to describe more accurately the stiffness behavior of most material continua \cite{scalerandi2016power}, whereas for various springs and composite meta-materials quadratic and cubic stiffness terms are appropriate \cite{hu2017nonlinear, tang2016using}. Since damping --and in particular hysteretic damping--crucially affects the dynamical behavior of these nonlinear models, we have decided to make this the focus of the present study, since research in this direction has been rather limited.
Thus, in this paper we discuss the numerical stability of two highly influential (and by no means simple) mathematical models for hysteretic damping in the presence of nonlinearity: the complex stiffness model first introduced by Bishop \cite{bishop1955treatment} and the internal friction model first introduced by Reid \cite{reid1956free}.

Many energy dissipation mechanisms arising in mechanical systems are characterized by nonlinearities of various types. In many studies, hysteretic damping is considered as the primary source of nonlinearity, even though other  forms of inherently nonlinear terms may also be present that can seriously affect the dynamics. In contrast with models that include nonlinear stiffness, systems whose only nonlinearity enters through damping are often called ``quasilinear'' \cite{elliott2015nonlinear} and are considered accurate approximations of the underlying dynamics, due to the absence of jump or bifurcation effects, which do not frequently occur in practice.  

In the absence of external forcing, if stiffness nonlinearities are small, their effect on the oscillator's behavior is expected to be minimal. However, in the case of periodic driving, the nonlinear terms can seriously affect the amplitude of the periodic solution to which the oscillations are attracted. In fact, under suitable parameter values, they may render the solution unstable, and steer the motion in a chaotic region, from which it may even escape to infinity, depending on the form of the potential (see e.g. \cite{strogatz2018nonlinear}).

Thus, one of the main objectives of this work, from a mechanical engineering point of view, is to ensure that hysteretically damped oscillators respond predictably under periodic forcing, exhibiting oscillations whose amplitude, frequency and phase can be efficiently controlled. This is especially significant, as we shall demonstrate, in the highly realistic case that such oscillators operate in the presence of weakly nonlinear stiffness terms.

In a seminal paper \cite{bishop1955treatment} R. E. D. Bishop introduced a complex hysteretic damping term for linear springs, as an improvement over the usual viscous damping models in which dissipation is simply proportional to the velocity of the oscillator $\dot{x}(t)$. In Bishop's hysteretic damping model the total (linear) force is {\it not in phase} with $x(t)$ and this accounts for dissipative effects due to hysteresis losses in the material of the spring. Thus both stiffness and damping are governed by properties of the material and are included in a single linear complex term in the equation of motion. 

In this paper, we are interested in adding a nonlinear stiffness term and driving such an oscillator by an external periodic force of the form $Fexp(i\omega t)$ to examine the periodic solution of period $T=2\pi/\omega$ to which the motion is attracted, since transients will be damped out in time. As Bishop demonstrated, in his linear stiffness model, the mathematical analysis is straightforward and can be easily carried out, not only in the one oscillator case \cite{bishop1955treatment}, but also for an an arbitrary number of $N$ coupled masses and springs \cite{bishop1956general}.

Since their introduction, hysteretic damping models have been used in a wide variety of applications, including seismic behavior \cite{zhu2007seismic}, composite beam modeling \cite{banks1991damping}, rotor dynamics \cite{montagnier2007dynamic, genta2004persistent} and material modeling \cite{gandhi1999characterization}, mainly because, contrary to the commonly used (fractional) Kelvin--Voigt based models \cite{lewandowski2010identification}, they allow for the energy loss per cycle to be  \textit{independent of the deformation frequency} \cite{inaudi1995linear}.

Our aim in this work is to study Bishop's hysteretic damping linear model in the presence of small stiffness nonlinearities, which typically occur in a wide variety of coupled oscillator systems of mechanical engineering \cite{strogatz2018nonlinear, Vakakisetal2008}. We thus begin in Section 2 by reviewing Bishop's linear oscillator model, pointing out certain spurious errors that arise in the numerical computation of the exact solution, already in the unforced case, at times when the oscillations become too small. We show that these errors persist even when we include periodic forcing on the right side of the equation and explain how they can be distinguished from the divergence of solutions due to dynamical effects.

Next, we proceed in Section 3 to solve the periodically driven single oscillator, in the presence of a small quadratic stiffness nonlinearity, by developing a perturbation expansion in powers of a small parameter $\epsilon$. Since part of the solution damps out due to dissipation, we focus on the periodic attractor to which the motion converges as $t\longrightarrow \infty$. This has the advantage of avoiding numerical errors produced by small oscillations, and allows us to verify computationally the validity of our analytical results, given as a convergent Fourier expansion of the $T$--periodic attractor, where $T$ is the period of the forcing term.

Studying in detail the behavior of the amplitude of the real part of our periodic solutions for different values of the parameters, including $\epsilon>0$, we find that it increases significantly in magnitude and complexity with growing $\epsilon$. Furthermore, plotting the amplitude response of our solutions as a function of frequency, we observe, besides the high peak at the harmonic frequency $\omega_1$ of the unforced oscillator, secondary peaks at subharmonics, $\omega_1 /2, \omega_1 /3,...$, which also grow with increasing $\epsilon$. We point out, however, that, as $\epsilon$ increases, the periodic solutions of the quadratic model eventually become unstable and lead to oscillations that escape to infinity through the homoclinic tangle associated with the saddle point of the potential.

Next, in Section 4, we turn to the case of a {\it symmetric} nonlinearity with positive quartic term in the potential and cubic nonlinearity in the equation of motion, for which saddles are absent and escape of solutions to infinity is precluded. We again obtain analytically the periodic solutions of the problem, and show that they become unstable as $\epsilon$ is increased, leading to chaotic motion very similar to what is observed in a classical Duffing's equation \cite{strogatz2018nonlinear}. This leads to oscillation amplitudes of great complexity in phase space projections, as well as amplitude and phase response curves with very similar features as in the case of quadratic nonlinearity. 

Nevertheless, as we demonstrate in this paper, in both nonlinear extensions, the solution of the Bishop model is never numerically stable for arbitrarily long times! Indeed, even in cases where our periodic solutions should be dynamically stable, they are {\it not} true attractors, since the numerical computation possesses inherent errors which will eventually grow exponentially and lead away from the desired periodic solution!

To address this problem, we devote Section 5 of our paper to a preliminary study of a different hysteretic model introduced by Reid \cite{reid1956free}, which appears to be free from the above numerical instabilities. Its linear version is expressed in terms of a differential equation with purely real (or purely imaginary) solutions and is hence more ``physical'' than Bishop's model. As a result, no numerical errors appear during integration and stable periodic solutions are true attractors. More importantly, the solutions of its nonlinear version exhibit a fascinating property of multistability, i.e. coexisting periodic attractors of different periods with intricate basins of attraction, whose rich structure needs to be further studied in the future.

Finally, in Section 6, we summarize our conclusions and describe ongoing work on the remarkable and potentially highly relevant dynamical features of nonlinear mechanical systems with hysteretic damping. The expected outcome of this work is to provide concrete and implementable criteria that can guarantee the safe and efficient operation, not only of one, but many hysteretically damped oscillators connected in a ``chain'' through nearest neighbor interactions.

Although the models discussed here are limited to the single-degree-of-freedom case, while practical engineering problems are usually delegated to multi-degree-of-freedom and eventually finite element analyses, most steady-state solvers for the latter are founded on these same models (especially Bishop’s model) and are therefore affected critically by the errors discussed explicitly for the first time in this paper. It is notable that, although the causality of Bishop’s model has been often been brought to question \cite{crandall1970role}, its numerical accuracy has barely received the same scrutiny. As such, we consider the findings presented in this paper to be significant.
At the same time, Reid’s model is one of the few truly practical alternatives to Bishop’s model among a plethora of other models, which use e.g. internal variables (Bouc-Wen) \cite{ikhouane2007dynamic}, fractional derivatives \cite{shen2020fractional}, Hilbert-transforms (Biot’s model) \cite{mastroddi2019time}, and even operators and neural networks \cite{dang2005neural}. These are mathematically interesting but practically too complex and sometimes artificial for general engineering practice. Reid’s model, on the other hand, has hardly received the attention it merits and several decades after its first introduction its uses remain limited to linear material models \cite{liu2006reid, kang2019steady, spitas2009continuous}. The present paper is perhaps the first to explore the use of this model for engineering materials with cubic stiffness nonlinearity to uncover significant behaviors that were previously unknown.

\section{A linear oscillator with hysteretic damping}

The analysis of one degree--of--freedom linear oscillator under periodic forcing in the presence of ``hysteretic damping'' \cite{bishop1955treatment} starts with the following complex differential equation:
\begin{equation}\label{bishop1}
M\ddot{x} + \dfrac{h}{\omega}\, \dot{x} + k x = F \exp{i\omega t},
\end{equation}
where $x$ denotes the particle displacement from equilibrium.

Note that when the coefficient of hysteretic damping $h$ is small compared to the stiffness $k$ the results approximate very closely those obtained for viscous damping. In particular, with $x = R\exp{i\omega t}$, we may write $h\dot{x}/\omega\rightarrow ihx$, and end up with the equation
\begin{equation}\label{bishop1b}
M\ddot{x} + (k + i h)\, x = F \exp{i\omega t},
\end{equation}
In this case, the complex stiffness of a hysteretic damper is given by $k + ih$ so that, for a spring-and-damper system, we have a (complex restoring) force $f=(k  + ih) x = k\,(1+i\mu)x$, with $\mu=h/k$, thus allowing for the desired frequency independent response.

\subsection{Free vibration of the linear oscillator}
Let us first analyze the free vibration of the linear hysteretic model, setting $F=0$ in Eqn. (\ref{bishop1b}) with $M=1$:
\begin{equation}\label{bishop2}
\ddot{x} + (k + i h)\, x = \ddot{x} + k\,(1 + i \mu)\, x = 0.
\end{equation}
We define the natural frequency of the oscillator by $\omega_1  = \sqrt{k}$, whence using the representation $x=u+iv$ we write the solution of Eqn. (\ref{bishop2}) in the form:
\begin{equation}\label{bishop2f}
x = \alpha \, \exp{-\omega_1bt + i \left(\omega_1ct+\theta\right)},
\end{equation}	
with real and imaginary parts
\begin{align}\label{bishop2g}
u &= \alpha\,\exp{-\omega_1bt}\cos{\left(\omega_1ct+\theta\right)}, \\
v &= \alpha\,\exp{-\omega_1bt}\sin{\left(\omega_1ct+\theta\right)},
\end{align}
respectively, where 
\begin{equation}\label{bishop2fa}
b = \left(\frac{\sqrt{1+\mu^2}-1}{2}\right)^{\frac{1}{2}},\,\,\, c = \left(\frac{\sqrt{1+\mu^2}+1}{2}\right)^{\frac{1}{2}}.
\end{equation}
This can be shown by directly differentiating the solution (\ref{bishop2f}) and substituting in Eqn. (\ref{bishop2}) to verify the solution as given in Eqs. (\ref{bishop2f}) and (\ref{bishop2fa}).

Setting now $t=0$ in Eqn. (\ref{bishop2f}) and its derivative we obtain the initial displacement and velocity as:
\begin{equation}\label{bishop2m}
x(0) = \alpha \exp{i\,\theta},\,\,\,\ \dot{x}(0) = \alpha\,\omega_1\,\left(-b + i \,c\right) \, \exp{i\theta}
\end{equation}
where the (real) parameters $\alpha$ and $\theta$ are arbitrarily chosen.

\begin{figure}[t]
	\centering
	\includegraphics[width = 0.5\textwidth]{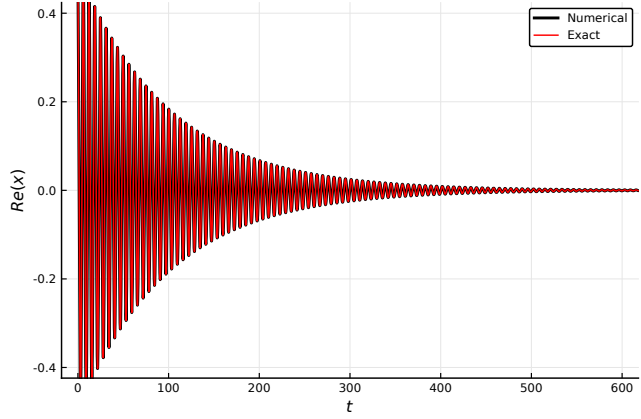}
	\caption{Using the parameters $k=1$, $h=0.02$ and solving Eqn. (\ref{bishop2}) for initial conditions with $\alpha = 0.5, \theta=0.3$ in Eqs. (\ref{bishop2m}) we find that numerical and exact solutions of the linear oscillator coincide up to time t=1000, in this example. For longer times, however, numerical errors are seen to grow exponentially.}	\label{fig2}
\end{figure}

It is instructive to solve Eqn. (\ref{bishop2}) numerically to see whether we can reproduce the analytical results for sufficiently long times. Thus, choosing initial conditions with $\alpha = 0.5, \theta=0.3$ in Eqn. (\ref{bishop2m}) and parameters $k=1$, $h=0.02$, we present in Fig. \ref{fig2} time plots of the displacement, superimposed with the exact solution for the real and imaginary parts of the displacement and verify that the numerical and analytical solutions practically coincide and appear to converge to zero, at least up to times $t\approx 1000$ in this example. 

However, the situation is quite different, if we repeat these calculations for the same initial conditions, $k=1$ and larger damping $h$, which causes the solutions to converge to zero faster. In fact, for $h=0.5$, the numerical and exact solutions are indistinguishable until a little before $t=100$, where numerical errors suddenly begin to grow. In fact, when we repeated the experiment of Fig. \ref{fig2} for times longer than t=1000, we also witnessed eventually a spurious ``blow-up'' of the numerical solutions!

We thus arrive at an important observation: While the numerical solution practically coincides with the analytical one up to times when they are both almost equal to zero, numerical instabilities arise, which lead to a spurious exponential blow-up. The critical time $t_d$ at which these errors begin to grow crucially depends on both the precision level of the computations and the specific choice of the parameter $\mu$. Interestingly enough, the dependence of $t_d$ on $\mu$ follows a power-law, which further supports our claim concerning the spurious nature of these errors.

The same type of errors occur during numerical simulations of the periodically forced linear problem as well as the nonlinear problems we will discuss in later sections. Thus, we henceforth adopt the following strategy when solving these equations numerically: When faced with diverging solutions starting at some approximate value $t=t_d$, we increase the accuracy of our computations and check whether $t_d$ depends on the specified tolerance. If that is the case, we trust the numerical results up the time of blow-up and conclude that the observed divergence is not a consequence of the true dynamical properties of the motion.

\subsection{Forced vibration}
Let us now briefly describe the solutions of the associated inhomogeneous problem adding a complex periodic force on the right hand side, as described in Eqn. (\ref{bishop1b}), where setting again $M=1$ we obtain
\begin{equation}\label{bishop1b2}
\ddot{x} + (k + i h)\, x = F \exp{i\omega t},
\end{equation}
with $F=f+i\,g$.

As is well--known, the solution of an inhomogeneous linear problem can be written the form $x = x_h + y$, with $x_h$ the solution of the associated homogeneous problem and $y$ a particular solution. Thus, substituting in Eqn. (\ref{bishop1b2}), we find:
\begin{equation}
\ddot{x_h} + \ddot{y} + \omega_1^2\,(1 + i \mu)\, \left(x_h + y\right) = F \exp{i\omega t},
\end{equation}
where, with $y=Bexp(i\omega t)$, we obtain
\begin{equation}
\ddot{y} + \omega_1^2\,(1 + i \mu)\, y = F \exp{i\omega t},\,\,\,\ B = \dfrac{F}{-\omega^2 + \omega_1^2\left(1+i\,\mu\right)}.
\end{equation}

So, we have the solution
\begin{equation}\label{insol}
x(t) = \alpha \, \exp{-\omega_1bt + i \left(\omega_1ct+\theta\right)} + B \exp{i\omega t}.
\end{equation}

In order to identify the initial conditions, we set $t = 0$ and obtain
\begin{equation}\label{incond}
x(0) = x_h(0) + B \,\,\, \text{and}\,\,\,\dot{x}(0) = \dot{x}_h(0) + i\omega B,
\end{equation}
for the initial displacement and velocity, respectively.

Plotting the magnification factor $n = \left(\left(1-\frac{\omega^2}{\omega_1^2}\right)^2+\mu^2\right)^{-0.5}$ with respect to the frequency ratio $\omega/\omega_1$ we find precisely the same results as obtained in \cite{bishop1955treatment}.

\begin{figure}[t]
	\centering
	\includegraphics[width = 0.65\textwidth]{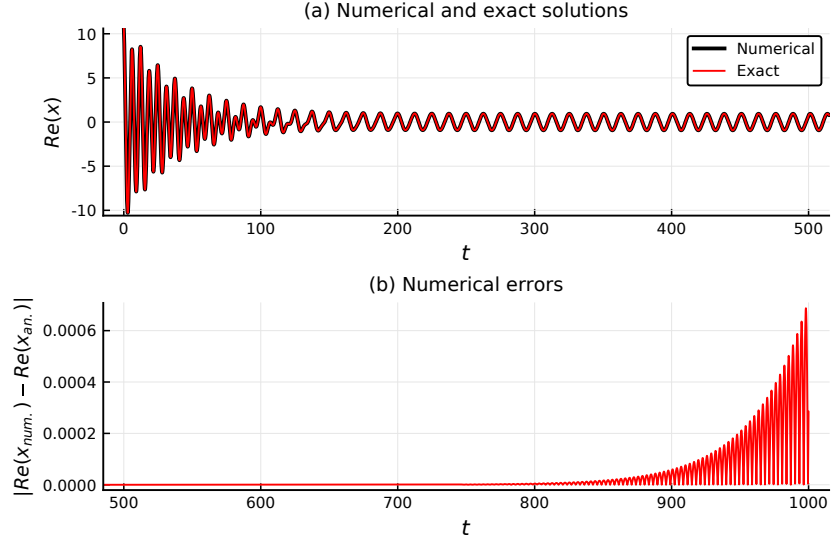}
	\caption{(a) Numerical and exact solutions of the inhomogeneous problem, together with (b) the corresponding errors. In his case, the correspondence is exact up to time $t\approx 800$, when the errors begin to grow. The parameters are set here to $(k,h,f,g,\omega) = (1,0.05,0.5,0.5,0.5)$, while the initial conditions correspond to $(\alpha,\theta) = (10.5,0.3)$}.	\label{fig4}
\end{figure}

Setting the parameters $(k,h,f,g,\omega) = (1,0.05,0.5,0.5,0.5)$ and for initial conditions associated with the initial parameters $(\alpha,\theta) = (10.5,0.3)$, we plot the results in Fig. \ref{fig4}, superimposing numerical and exact solutions. After a transient time, the trajectory (of the real part of the displacement) oscillates periodically with frequency equal to the driving frequency $\omega$.

Note, however, in Fig. \ref{fig4}, just as in the unforced case, that our numerical solution reproduces accurately the analytical results up to a time interval $t_d\approx 800$, as the damped part of the solution produces numerical errors, which grow exponentially (see Fig. \ref{fig4}(b)). Thus, we verify here also the presence of the same type of unavoidable errors as in the unforced problem, whose occurrence depends on the specified accuracy of the computations.  

\section{The case of nonlinear stiffness}

We now come to the main part of our study, where we examine the effect of adding a nonlinear stiffness term to our model to explore its influence on the dynamics. As is well--known, realistic springs are in general not linear, and may include nonlinearities that depend on the displacement. Thus, we will assume at first that these nonlinearities are quadratic and later examine also the case of {\it symmetric} springs, where the lowest order nonlinearities are cubic in the displacement. More specifically, our nonlinear hysteretic damping model has the form
\begin{equation}\label{nonlin1}
\ddot{x} + (k + i h)\, x + \epsilon \, x^2 = \ddot{x} + \omega_1^2 \, (1 + i \mu)\, x + \epsilon \, x^2 = F \exp{i\omega t},
\end{equation}
with $F$ real.

Furthermore, since the damped part of the solution vanishes exponentially after relatively small time intervals, we focus on the periodic attractor to which the motion eventually converges and study the effect of nonlinearity on its amplitude and frequency response. This permits us to avoid the spurious numerical errors caused by the damping and concentrate on the ultimate periodic solution, whose form can be obtained analytically as a Fourier series that can be reliably tested against the numerical solution. 

To this end, we begin by applying perturbation theory \cite{wiggins2003introduction} to derive the first few terms of the solution in powers of a small parameter $\epsilon$, and then develop the full Fourier series that represents the desired periodic solutions. Next, we demonstrate the convergence of this series and compare its accuracy against numerical results.

\subsection{Analytical results}

We seek an approximate solution of Eqn. (\ref{nonlin1}) as an asymptotic series in powers of $\epsilon$ of the following form:
\begin{equation}\label{app_solution}
x(t) = x_0(t) + \epsilon \, x_1(t) + \epsilon^2 \, x_2(t) + \ldots
\end{equation}

Substituting this expression in Eqn. (\ref{nonlin1}), we equate terms of like powers of $\epsilon$ and obtain, as usual, an infinite sequence of linear inhomogeneous equations for the $x_j$, with identical homogeneous part, as follows:
\begin{align}
\epsilon^0: \,\,\, \ddot{x_0} + \omega_1^2 \, (1 + i \mu)\, x_0 &= F \exp{i\omega t}\label{e0} \\
\epsilon^1: \,\,\, \ddot{x_1} + \omega_1^2 \, (1 + i \mu)\, x_1 &=  -x_0^2\label{e1} \\
\epsilon^2: \,\,\, \ddot{x_2} + \omega_1^2 \, (1 + i \mu)\, x_2 &= -2 x_0 x_1 \label{e2} \\
\epsilon^3: \,\,\, \ddot{x_3} + \omega_1^2 \, (1 + i \mu)\, x_3 &= -x_1^2 - 2 x_0 x_2,\label{e3}
\end{align}
etc. From the lowest order terms we recover the inhomogeneous second order linear equation, given by Eqn. (\ref{insol}), whose solution we write here in the form:
\begin{equation} \label{x0}
x_0(t) = \alpha \, \exp{g(t)} + B_0 \exp{i\omega t},\,\,\,B_0 = \dfrac{F}{-\omega^2 + \omega_1^2\left(1+i\,\mu\right)}
\end{equation}	
where $g(t) = -\omega_1bt + i \left(\omega_1ct+\theta\right)$.

For the term of order $\epsilon^1$ we have:
\begin{align}
\ddot{x_1} + \omega_1^2 \, (1 + i \mu)\, x_1 = &- \alpha^2\, \exp{2 g(t)} - \label{x1.1}\\
 &2\alpha B_0 \exp{g(t)+i\omega t} - B_0^2 \exp{2 i\omega t}, \nonumber
\end{align}
whose solution is of the form:
\begin{equation}\label{x1.2}
x_1(t) = \lambda_1 \exp{2 g(t)} + \lambda_2 \exp{g(t)+i\omega t} + B_1 \exp{2i\omega t},
\end{equation} 
Substituting (\ref{x1.2}) in (\ref{x1.1}) we derive the following expressions for $\lambda_1,\lambda_2$: 
\begin{align}
\lambda_1 &= \dfrac{-a^2}{4\left(-\omega_1 b + i\omega_1c\right)^2 + \omega_1^2(1+i\mu)} \\
\lambda_2 &= \dfrac{2\alpha F}{\left(-\omega^2 + \omega_1^2(1+i\mu)\right) \left(\left(-\omega_1 b + i\omega_1 c + i\omega\right)^2+\omega_1^2(1+i\mu)\right)}\label{x1l1}
\end{align}
while for $B_1$ we obtain:
\begin{equation}
B_1 = \dfrac{-F^2}{\left(\Omega_1^2 - \omega^2 \right)^2 \left(\Omega_1^2 - 4\omega^2 \right)} \label{x1l3},
\end{equation}
using the notation $\Omega_1^2 \equiv \omega_1^2(1+i\mu)$.

For the $\epsilon^2$ term we have:
\begin{align}
&\ddot{x_2} + \omega_1^2 \, (1 + i \mu)\, x_2 =  -2 x_0 x_1 \nonumber\\ 
= &-2 \left(\dfrac{F\exp{i\,\omega\,t}}{-\omega^2 + \omega_1^2\left(1+i\,\mu\right)}\right) \left(\dfrac{-F^2 \exp{2i\,\omega\,t}}{\left(\Omega_1^2 - \omega^2 \right) \left(\Omega_1^2 - 4\omega^2 \right)}\right) \nonumber \\
&= B_2 \exp{3i\omega t},\label{x2.1}
\end{align}
where we kept only the $\exp{i\omega t}$--part of the $x_0$ solution and the $\exp{2 i\omega t}$--part of the $x_1$ solution. The reason for this, is that from now on we will neglect the damped part and focus on the periodic attractor, for which we will construct its Fourier series representation. To this end, we first identify the coefficient $B_2$ of the term $\exp{3i\omega t}$ in the above calculation as:
\begin{equation}
B_2 = \dfrac{2 F^3}{\left(\Omega_1^2 - \omega^2 \right)^3 \left(\Omega_1^2 - 4\omega^2 \right) \left(\Omega_1^2 - 9\omega^2 \right)}. \label{x2b2}
\end{equation}

Proceeding finally to the terms of order $\epsilon^3$, we find from (\ref{e3}):
\begin{equation}
\ddot{x_3} + \omega_1^2 \, (1 + i \mu)\, x_3 = (-B_1^2-2 B_0\,B_2) \exp{4i\omega t},\label{x3.1}
\end{equation}
and obtain a solution of the form $x_3 = B_3\exp{4i\omega t}$, which when substituted in Eqn. (\ref{x3.1}) yields:
\begin{equation}
B_3 = \dfrac{F^4}{\Omega_1^2-16\omega^2} \cdot \dfrac{\left(\Omega_1^2-9\omega^2\right) + 4\left(\Omega_1^2-4\omega^2\right)}{\left(\Omega_1^2-\omega^2\right)^4 \left(\Omega_1^2-4\omega^2\right)^2 \left(\Omega_1^2-9\omega^2\right)}. \label{x3.2}
\end{equation}
In a similar manner, we proceed to the next higher order terms and evaluate the coefficients $B_4$ and $B_5$ for the associated solutions $x_4 = B_4\exp{5 i\omega t}$ and $x_5 = B_5\exp{6 i\omega t}$, and so on. Thus, we write the complete Fourier series expansion of the periodic solution $\hat{x}(t)$ in the form:
\begin{equation}
\hat{x}(t) = \sum_{k=1}^{\infty}\epsilon^{k-1}B_{k-1}\exp{k i \omega t},\label{Foursol}
\end{equation}
whose coefficients are recursively given by the expression:
\begin{equation}\label{rec_eq}
B_k = \dfrac{-1}{\Omega_1^2-\left[(k+1)\omega\right]^2}\left(B^2_{\nu-1}\delta_{k,2\nu-1} + 2\sum_{\substack{i+j=k-1\\
		i>j}}B_iB_j\right),
\end{equation}
for $k\geq1, \nu\geq2$, where $B_0 = \dfrac{F}{\Omega_1^2-\omega^2}$ and $\delta(n,m)$ is the Kronecker delta that equals 1 for $n=m$ and 0 for $n\neq m$.

\begin{figure}[t]
	\centering
	\includegraphics[width = 0.65\textwidth]{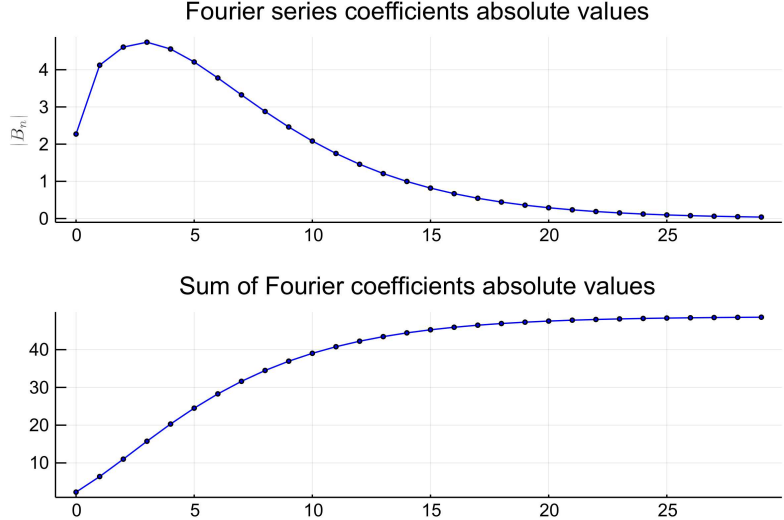}
	\caption{Illustration of the convergence of the Fourier series expansion of the periodic solution for the nonlinear problem. Parameters used: $\left(\mu,f,g,\omega,\epsilon\right) = \left(0.05, 1, 0, 0.75, 0.1\right)$.}	\label{fig6a}
\end{figure}

Fig. \ref{fig6a} illustrates the absolute convergence of the obtained Fourier series, for the choice of parameters $\left(\mu,f,g,\omega,\epsilon\right) = \left(0.05, 1, 0, 0.75, 0.1\right)$. This is supported by the observation that the $|B_n|$ coefficients in the top graph of Fig. \ref{fig6a}, decay like $n^{-3}$ for large $n$.

\subsection{Stability of solutions and amplitude and phase response curves}

It is instructive to examine the accuracy of our analytical results for the nonlinear system by comparing the solutions obtained numerically and analytically, using the recursive relation (\ref{rec_eq}) for the Fourier coefficients. First, we examine the validity of our series expansion by keeping 150 terms in (\ref{Foursol}) and evaluating the real and imaginary parts of the quantity:
\begin{equation}
Q = \hat{\ddot{x}} + \omega_1^2 \, (1 + i \mu)\, \hat{x} + \epsilon \, \hat{x}^2 - F \exp{i\omega t},
\end{equation}
which is expected to be \textit{close to zero}, up to machine precision. The results are quite encouraging as we find, for $\left(\mu,f,g,\omega\right) = \left(0.05, 1, 0, 0.75\right)$ and $\epsilon=0.01$ that the above quantity $Q$ oscillates in the range $\left[-2\times10^{-15},2\times10^{-15}\right]$, at least up to $t=1000$ units.

\begin{figure}[t]
	\centering
	\includegraphics[width = 0.65\textwidth]{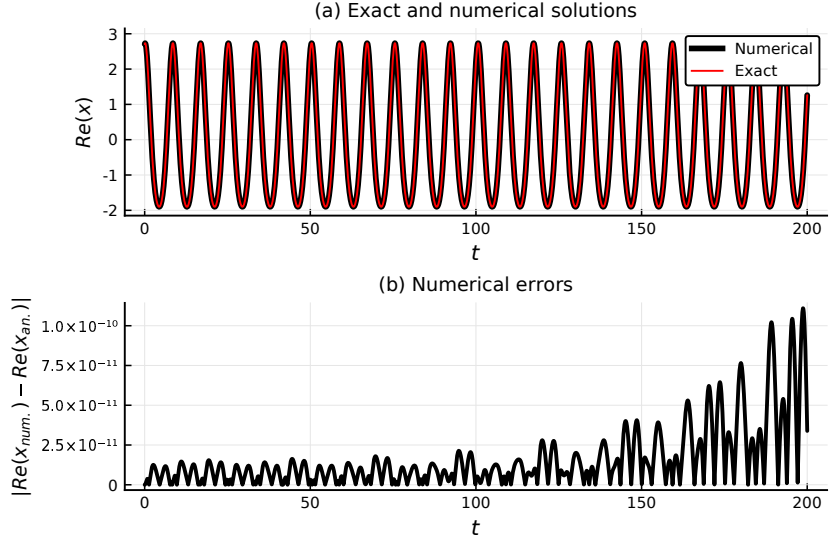}
	\caption{Top: Numerical (black) and analytical (red) trajectory of the nonlinear system of Eqn. (\ref{nonlin1}) for $\epsilon=0.01$. Bottom: Time plot of the ``error'', calculated as the distance between the numerical and the analytical solution, i.e. $\lvert \text{Re}(x_{\text{an}}) - \text{Re}(x_{\text{num}}) \rvert$. For the analytical solution we kept the first 150 terms of the Fourier expansion.}	\label{fig6}
\end{figure}

We now increase $\epsilon$ to 0.1 and compare in Fig. \ref{fig6} (top) the numerical solution (in black), with the analytical solution (in red), for the same parameters $(\mu,f,g,\omega)$ as above. The two orbits are indistinguishable, but only for a time $t_d\approx 200$ as we verify by plotting in Fig. \ref{fig6} (bottom) the absolute error between them, $\lvert \text{Re}(x_{\text{an}}) - \text{Re}(x_{\text{num}}) \rvert$.

Next, we consider different values of $\epsilon \in \left\{0, 0.1, 0.25\right\}$ and observe that with increasing $\epsilon$ the periodic solutions become more ``sharply peaked'', while their amplitude reaches higher values. This is very important in view of the following argument: Since the stiffness terms of our nonlinear equation are connected with the potential function  
\begin{equation}
V(x)=\frac{1}{2}x^2 + \frac{\epsilon}{3}x^3 ,
\end{equation}\label{quadpot}
as is well--known from the corresponding {\it unforced} classical oscillator, if the solution amplitudes exceed the location of the saddle point of the above potential $x_s=-1/\epsilon$, the motion will diverge to $-\infty$. In the periodically forced case, this already occurs at values $x(t)<x_s$, since the motion in the neighborhood of $x_s$ is chaotic and orbits may well be driven to divergence at much smaller amplitudes \cite{strogatz2018nonlinear}. This is because, in such cases, the size of the region of regular motion, where stable periodic solutions exist, becomes severely limited, as the magnitude of the forcing amplitude increases, for any value of $\epsilon$.

\begin{figure}[t]
	\centering
	\includegraphics[width = 0.65\textwidth]{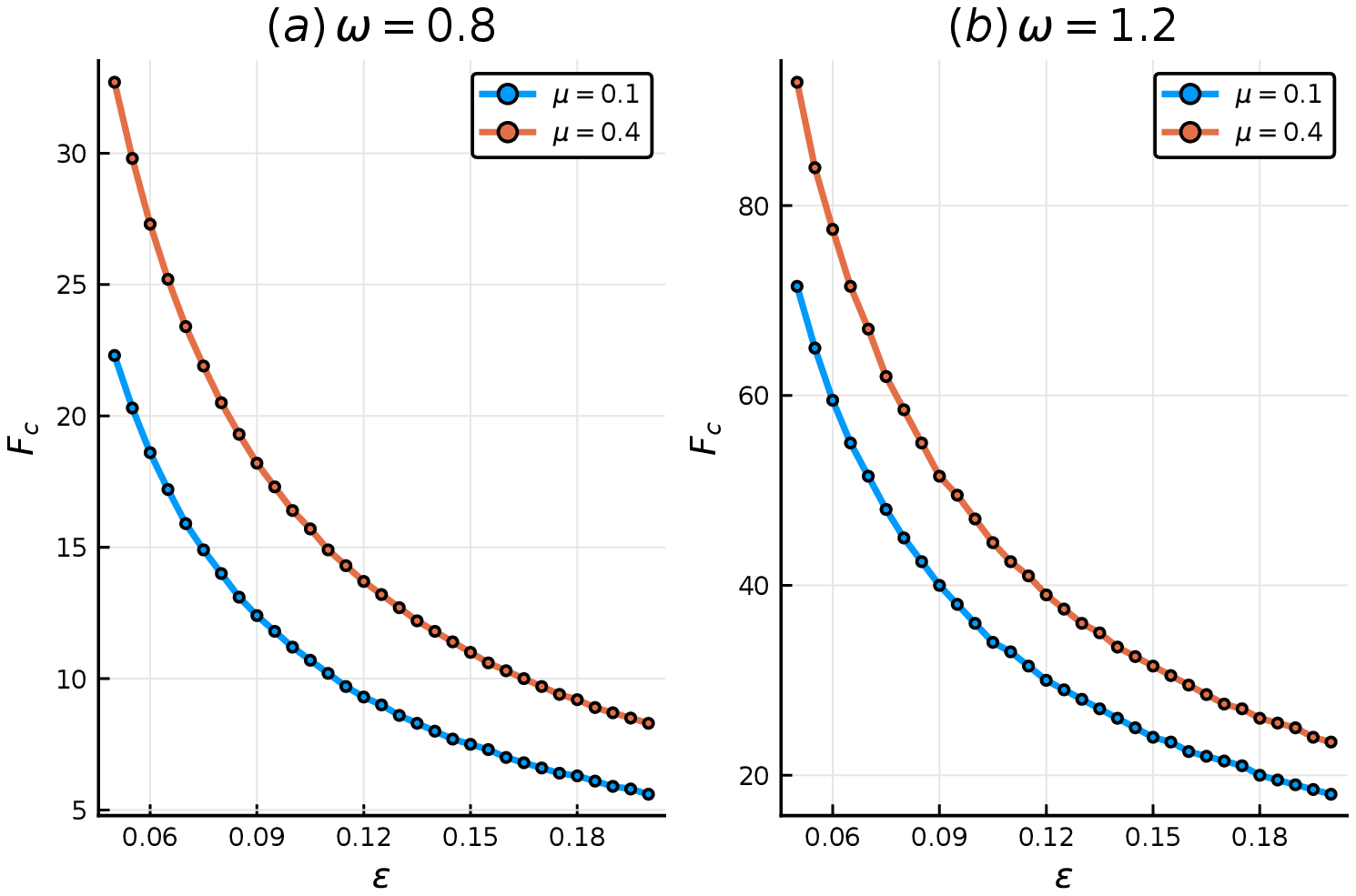}
	\caption{Critical forcing amplitudes $F_c$ with respect to the parameter $\epsilon$ for the system with quadratic nonlinearity, for frequencies  $\omega=0.8$ (left) and $\omega=1.2$ (right), letting $\mu \in \left\{0.1,0.4\right\}$.}	\label{fig13}
\end{figure}

It is, therefore, important to examine the critical forcing amplitude $F_c$ (leading to instability and escape of our solutions) in relation with $\epsilon$. As we see in Fig. \ref{fig13}, the increase of $\epsilon$ leads to lower critical amplitudes, thus making it ``easier'' for the orbits to escape. Moreover, we need to consider the effects of the forcing frequency $\omega$ and the measure of damping $\mu$ ($=h/k$). Specifically, higher forcing frequencies or higher damping measures both lead to higher critical amplitudes, thus making the oscillator more ``robust'' in terms of escaping.

\begin{figure}[t]
	\centering
	\includegraphics[width = 0.65\textwidth]{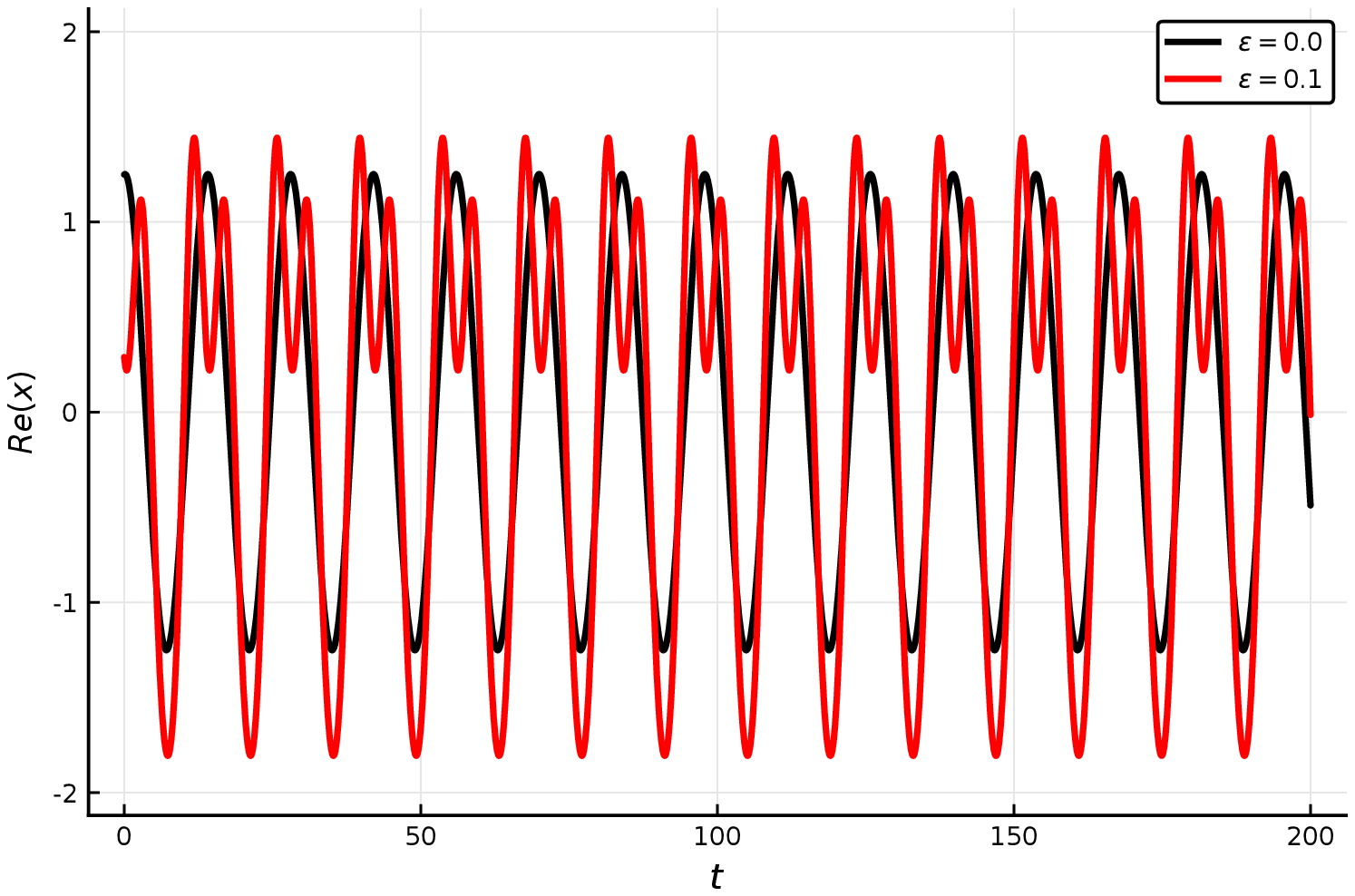}
	\caption{Trajectories of $\text{Re}(x)$ for different values of $\epsilon \in \left\{0 \text{(black)}, 0.1 \text{(red)}\right\}$, with $\left(\mu,f,g,\omega\right) = \left(0.05, 1, 0, 0.45\right)$.}	\label{fig7b2}
\end{figure}

As was already mentioned, the amplitude and morphology of the solutions depends significantly on the value of $\epsilon$, which may be crucial in certain mechanical applications. To further illustrate this important effect, we depict in Fig. \ref{fig7b2}, using $\left(\mu,f,g,\omega\right) = \left(0.05, 1, 0, 0.45\right)$, the linear ($\epsilon=0$) solution (black) together with the nonlinear one (red) at $\epsilon=0.1$. Notice that, while the period is the same, as expected, the nonlinear attractor exhibits distinctly larger amplitude and has a different shape when it reaches negative values. 

A related very important feature of our nonlinear system is the behavior of the amplitudes and phases of the periodic solution as functions of the driving frequency, since we now have more frequencies present in the solution due to the nonlinear term of the equation. In particular, we are interested in both amplitude and phase response curves and their dependence on the crucial parameters of the system. 

The amplitude response curve (magnification factor) $n$ and phase response curve $\eta$ are given by the expressions \cite{shabana2018theory,bishop1955treatment}:
\begin{equation}
n =\left(\left(1 - \frac{\omega^2}{\omega_1^2}\right)^2 + \mu^2 \right)^{-1/2},\,\,\ \eta = \tan[-1](\mu/\left(1 - \frac{\omega^2}{\omega_1^2}\right))
\end{equation} 
for the linear hysteretic damping model (i.e. $\epsilon=0$), and feature in its solution in the form
\begin{equation}
x = n \dfrac{F}{k}\exp{i (\omega t - \eta)},
\end{equation}
as functions of the frequency ratio $\omega/\omega_1$. Regarding the nonlinear model, however, since we do not have closed form expressions, we will compute the amplitude and phase response curves using the corresponding quantities obtained from our analytical solutions, $n=$ amplitude of oscillations/forcing amplitude and $\eta=$phase of oscillations$-$forcing phase.

\begin{figure}[t]
	\centering
	\includegraphics[width = 0.65\textwidth]{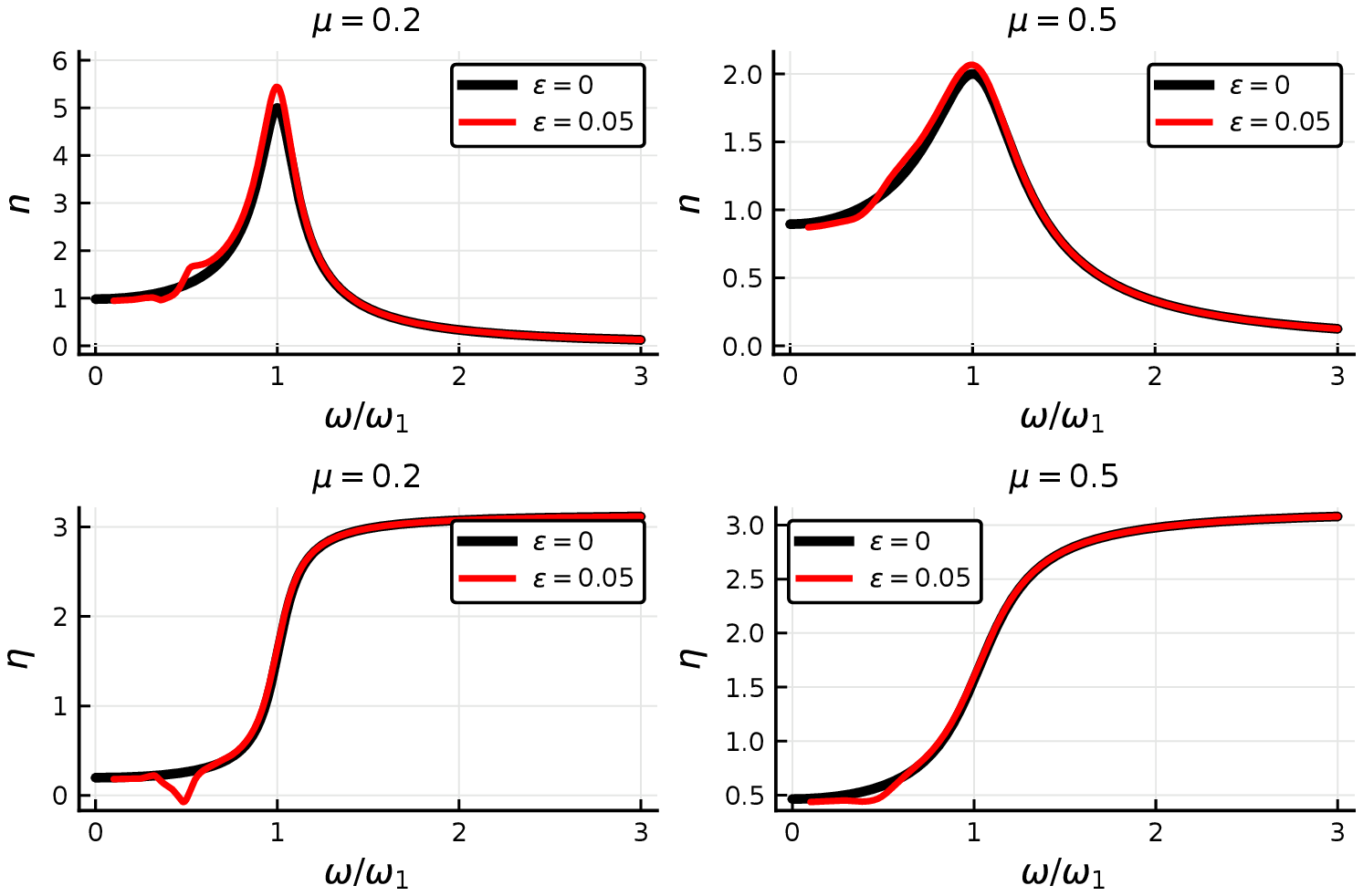}
	\caption{Amplitude (top) and phase (bottom) response curves of Eqn. (\ref{nonlin1}) for
		$\mu \in \left\{0.2, 0.5\right\}$ and $\epsilon \in \left\{0, 0.05\right\}$. With black lines we indicate the location of the frequencies $\omega_1/2$ and $\omega_1/3$.}	\label{fig8a}
\end{figure}

In Fig. \ref{fig8a} we plot the amplitude (top) and phase (bottom) response curves of Eqn. (\ref{nonlin1}) for $\mu \in \left\{0.2, 0.5\right\}$ and $\epsilon \in \left\{0, 0.05\right\}$.

\begin{figure}[t]
	\centering
	\includegraphics[width = 0.65\textwidth]{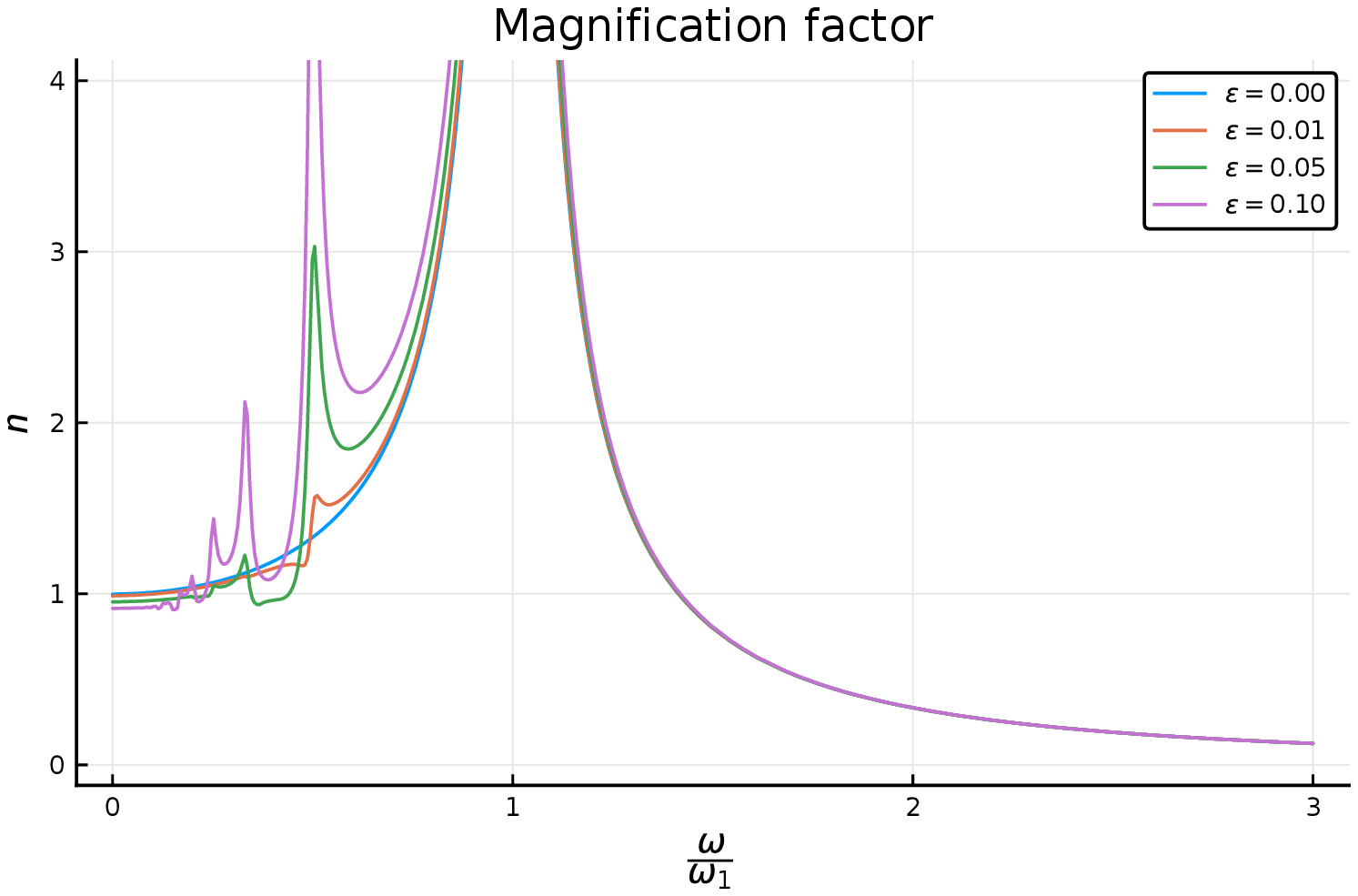}
	\caption{Amplitude response curves (magnification factor) of Eqn. (\ref{nonlin1}) for $\left(\mu,f,g\right) = \left(0.05, 1, 0\right)$ and $\epsilon \in \left\{0, 0.01,0.05,0.1\right\}$.}	\label{fig8b}
\end{figure}

In Fig. \ref{fig8b} we present amplitude response curves, keeping $\left(\mu,f,g\right) = \left(0.05, 1, 0\right)$ and varying $\omega$ for different values of  $\epsilon \in \left\{0, 0.01,0.05,0.1\right\}$. Notice the enhancement of the amplitude at the resonances $\omega \in\left\{\omega_1,\omega_1/2,\omega_1/3,\ldots\right\}$ as $\epsilon$ increases.

\section{Results for cubic nonlinearity}

After examining the case of quadratic nonlinearities, we now move on to {\it symmetric springs} where the nonlinearity is cubic in the displacement. More specifically, our nonlinear hysteretic damping model has the form
\begin{equation}\label{nonlin1b}
\ddot{x} + (k + i h)\, x + \epsilon \, x^3 = \ddot{x} + \omega_1^2 \, (1 + i \mu)\, x + \epsilon \, x^3 = F \exp{i\omega t}.
\end{equation}
As in the previous case, since the damped part of the solution vanishes exponentially after relatively small time intervals, we concentrate on the periodic attractor to which the motion eventually converges and study the effect of nonlinearity on its amplitude and frequency response. 

We begin again by applying perturbation theory to derive the first few terms of the solution in powers of a small parameter $\epsilon$, and develop the full Fourier series that represents the desired periodic solution. Next, we demonstrate the convergence of this series and compare its accuracy against numerical results.

\subsection{Analytical results}
We seek an approximate solution of Eqn. (\ref{nonlin1}) as an asymptotic series in powers of $\epsilon$ of the following form:
\begin{equation}\label{app_solution3}
x(t) = x_0(t) + \epsilon \, x_1(t) + \epsilon^2 \, x_2(t) + \ldots
\end{equation}
In what follows, we make use of the following expansion of a finite sum raised to the third power:
\begin{equation}\label{cub_exp}
\left(\sum_{i=1}^{N}\alpha_i\right)^3 = \sum_{i=1}^{N}\alpha_i^3 + 3\left(\sum_{i<j}^{N}\alpha_i \alpha_j^2 + \sum_{i<j}^{N}\alpha_i^2 \alpha_j\right) + 6\sum_{i<j<l}^{N}\alpha_i \alpha_j \alpha_l.
\end{equation}
Substituting this expression in Eqn. (\ref{nonlin1b}), we equate terms of like powers of $\epsilon$ and obtain, as usual, an infinite sequence of linear inhomogeneous equations for the $x_j$, with identical homogeneous part, as follows:
\begin{align}
\epsilon^0: \,\,\, \ddot{x_0} + \omega_1^2 \, (1 + i \mu)\, x_0 &= F \exp{i\omega t}\label{e03} \\
\epsilon^1: \,\,\, \ddot{x_1} + \omega_1^2 \, (1 + i \mu)\, x_1 &=  -x_0^3\label{e13} \\
\epsilon^2: \,\,\, \ddot{x_2} + \omega_1^2 \, (1 + i \mu)\, x_2 &= -3 x_0^2 x_1 \label{e23} \\
\epsilon^3: \,\,\, \ddot{x_3} + \omega_1^2 \, (1 + i \mu)\, x_3 &= -3 (x_0^2 x_2 + x_1^2 x_0)\label{e33}\\
\epsilon^4: \,\,\, \ddot{x_4} + \omega_1^2 \, (1 + i \mu)\, x_4 &= -\left(x_1^3 + 3 x_0^2 x_3 + 6 x_0 x_1 x_2\right)\label{e4}
\end{align}
etc. From the lowest order terms we recover the inhomogeneous second order linear equation, whose solution (neglecting the damping part) is written in the form:
\begin{equation} \label{x03}
x_0(t) =  B_0 \exp{i\,\omega\,t},
\end{equation}	
where $B_0 = \dfrac{F}{\omega_1^2\left(1+i\,\mu\right) -\omega^2} = \dfrac{F}{\Omega_1^2 -\omega^2}$. For the term of order $\epsilon^1$ we thus have:
\begin{equation}
\ddot{x_1} + \omega_1^2 \, (1 + i \mu)\, x_1 = - B_0^3 \exp{3 i\,\omega\,t} , \label{x1.13}
\end{equation}
whose solution is of the form:
\begin{equation}\label{x1.23}
x_1(t) = B_1 \exp{3 i\,\omega\,t},\,\,\ B_1 = -\dfrac{1}{\Omega_1^2 - (3\omega)^2} B_0^3,
\end{equation}
where we have used the notation $\Omega_1^2 \equiv \omega_1^2(1+i\mu)$.

For the $\epsilon^2$ term we get:
\begin{align}
&\ddot{x_2} + \omega_1^2 \, (1 + i \mu)\, x_2 =  -3 x_0^2 x_1 \nonumber \\
&= -3 B_0^2 B_1 \exp{5i\omega t}= B_2 \exp{5i\omega t},\label{x2.13}
\end{align}
from which we obtain $B_2$, while at order $\epsilon^3$ we also find $B_3$:
\begin{equation}
B_2 = \dfrac{-3B_0^2 B_1}{\Omega_1^2 - (5\omega)^2} ,\,\,\, B_3 = \dfrac{-3\left(B_0^2 B_2 + B_1^2 B_0\right)}{\Omega_1^2 - (7\omega)^2} . \label{x2b23}
\end{equation}

In a similar manner, we evaluate from the next higher order equations the coefficients $B_4$, $B_5$, etc.  write the complete Fourier series expansion of the periodic solution $\hat{x}(t)$ in the form:
\begin{equation}
\hat{x}(t) = \sum_{k=0}^{\infty}\epsilon^{k}B_{k}\exp{(2k+1) i \omega t},\label{Foursol3}
\end{equation}
whose coefficients are recursively given by the expression:
\begin{align}\label{rec_eq3}
B_k &= \dfrac{-1}{\Omega_1^2-\left[(2k+1)\omega\right]^2}\left(B^3_{\frac{k-1}{3}}\delta_{k,3\nu} + 3\sum_{\substack{i+2j=k-1\\
		i>j}}B_iB_j^2 + \right. \nonumber\\
	&\left. 3\sum_{\substack{2i+j=k-1\\
		i>j}}B_i^2 B_j + 6\sum_{\substack{i+j+l=k-1 \\
		i>j}}B_iB_jB_l\right),
\end{align}
for $k\geq1$, with $B_0 = \dfrac{F}{\Omega_1^2-\omega^2}$ and $\delta_{n,m}$=1 for $n=m$ and 0 for $n\neq m$.

\subsection{Numerical results}

\begin{figure}[t]
	\centering
	\includegraphics[width = 0.65\textwidth]{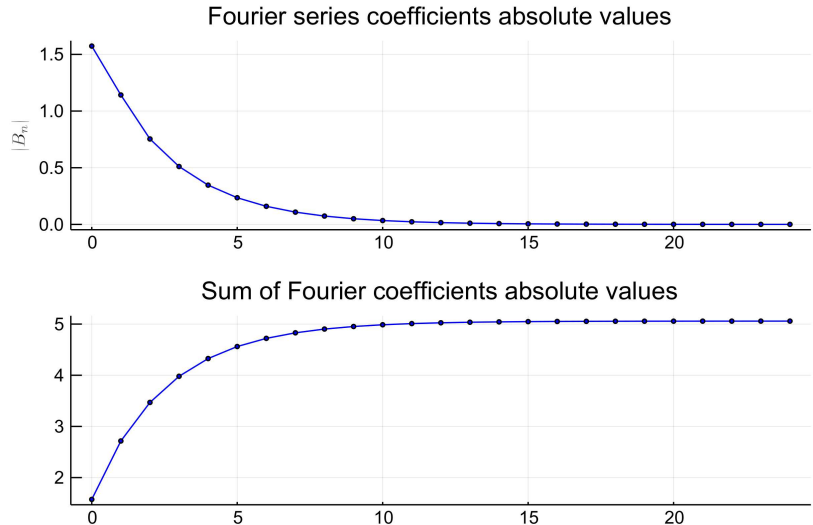}
	\caption{Illustration of the (absolute) convergence of the Fourier series expansion of the periodic solution for the nonlinear problem. Parameters used: $\left(\mu,f,g,\omega,\epsilon\right) = \left(0.05, 0.8, 0.1, 0.7, 0.1\right)$.}	\label{fig14}
\end{figure}

As in the case of the quadratic nonlinearity, it is possible to verify here also the absolute convergence of our Fourier series, by plotting the coefficients $|B_n|$ and their sum  in Fig. \ref{fig14} for the choice of parameters $\left(\mu,f,g,\omega,\epsilon\right) = \left(0.05, 0.8, 0.1, 0.7, 0.1\right)$.

\begin{figure}[t]
	\centering
	\includegraphics[width = 0.65\textwidth]{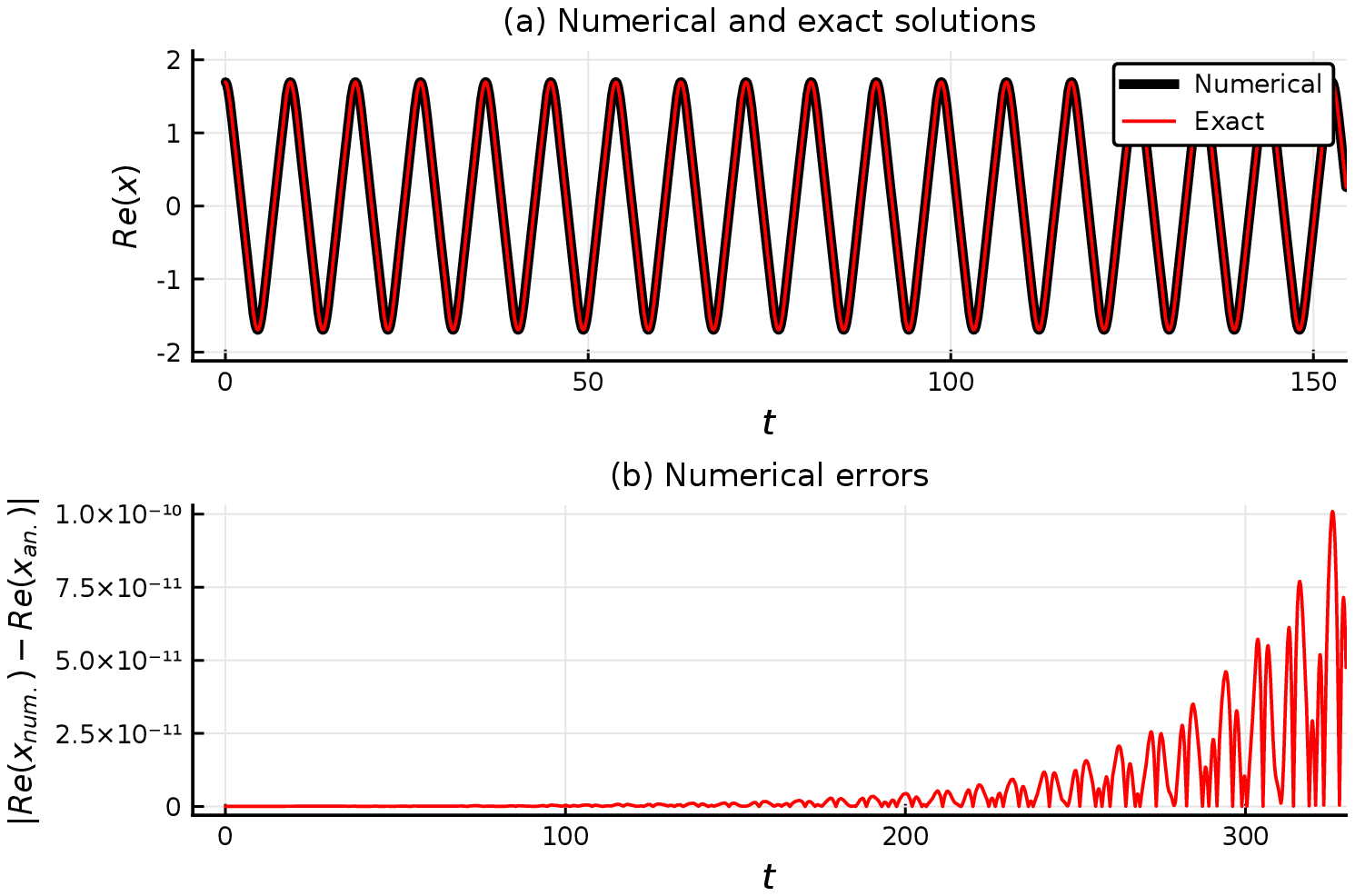}
	\caption{Top: Numerical (black) and analytical (red) trajectory of the nonlinear system of Eqn. (\ref{nonlin1}) for $\epsilon=0.1$. Bottom: Time plot of the error, calculated as the distance between the numerical and the analytical solution, i.e. $\lvert \text{Re}(x_{\text{an}}) - \text{Re}(x_{\text{num}}) \rvert$. For the analytical solution we kept the first 150 terms of the Fourier expansion.}	\label{fig16}
\end{figure}

Let us now examine the accuracy of our Fourier series, by comparing in Fig. \ref{fig16}(top) the numerical with the analytical solution, for the parameters$\left(\mu,f,g,\omega\right) = \left(0.05, 0.8, 0.1, 0.7\right)$ and $\epsilon = 0.1$. Initially, the two orbits are indistinguishable, as we can also also verify by the absolute error between them, defined $\lvert \text{Re}(x_{\text{an}}) - \text{Re}(x_{\text{num}}) \rvert$ in Fig. \ref{fig16} (bottom). However, as in the case of the quadratic nonlinearity, numerical instabilities also appear here at some value of $t=t_d$ that depends on the accuracy of the computations.

It is important to emphasize, however, that, unlike the model with quadratic nonlinearity, the equation with a cubic nonlinearity (\ref{nonlin1b}) is not expected to have solutions that diverge to infinity, since its equation is associated with the potential function
\begin{equation}
V(x)=\frac{1}{2}x^2 + \frac{\epsilon}{4}x^4 ,
\end{equation}\label{quadpot}
which has no saddle points, hence all motion is bounded. This implies that all solutions of this model must be dynamically stable and any kind of divergence, such as the one observed in Fig. \ref{fig16} must be of numerical origin. We have verified that this is indeed the case by checking that the time $t_d$ where our solutions begin to diverge depends on the selected precision of our computations.

\begin{figure}[t]
	\centering
	\includegraphics[width = 0.65\textwidth]{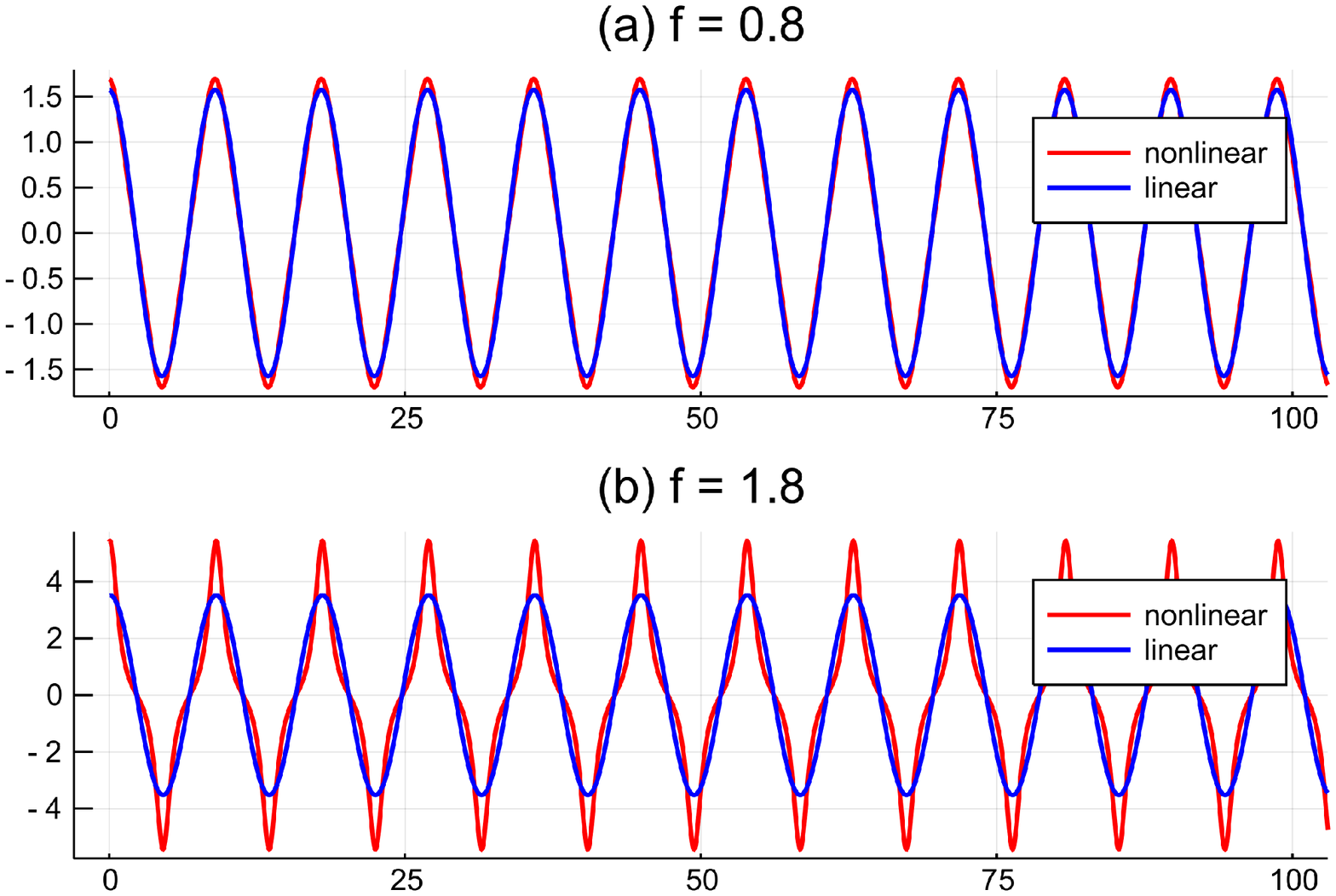}
	\caption{Time plots of real parts for the nonlinear ($\epsilon=0.1$) and the associated linear system ($\epsilon=0$), for different values of the (real part of the) driving amplitude $f$, namely $f=0.8$ (top) and $f=1.8$ (bottom). The rest of the parameters used $\left(\mu,g,\omega\right) = \left(0.05, 0.1, 0.7\right)$.}	\label{fig17}
\end{figure}

In order to study the differences between the linear and nonlinear system, we compare in Fig. \ref{fig17} the real parts of the associated solutions for $\epsilon = 0.1$ (red) and $\epsilon = 0$ (blue), using different values of the real part of the driving amplitude. As we can see, for small driving amplitude, the two solutions are nearly indistinguishable, except for a small amplitude difference. For $f=1.8$, however, keeping every other parameter fixed as above, we observe a distinct difference in the \textit{shape} of the periodic solution. 

\subsection{Comparison of the proposed models}

In this section we highlight the differences between of the 3 models, namely the linear and the ones with quadratic and cubic nonlinearity, in terms of the morphology of the periodic solutions in phase space as well as their amplitude and phase response curves.

\subsubsection{Morphology of the periodic solutions}

\begin{figure*}[t]
	\centering
	\begin{subfigure}[t]{0.5\textwidth}
		\centering
		\includegraphics[height=2in]{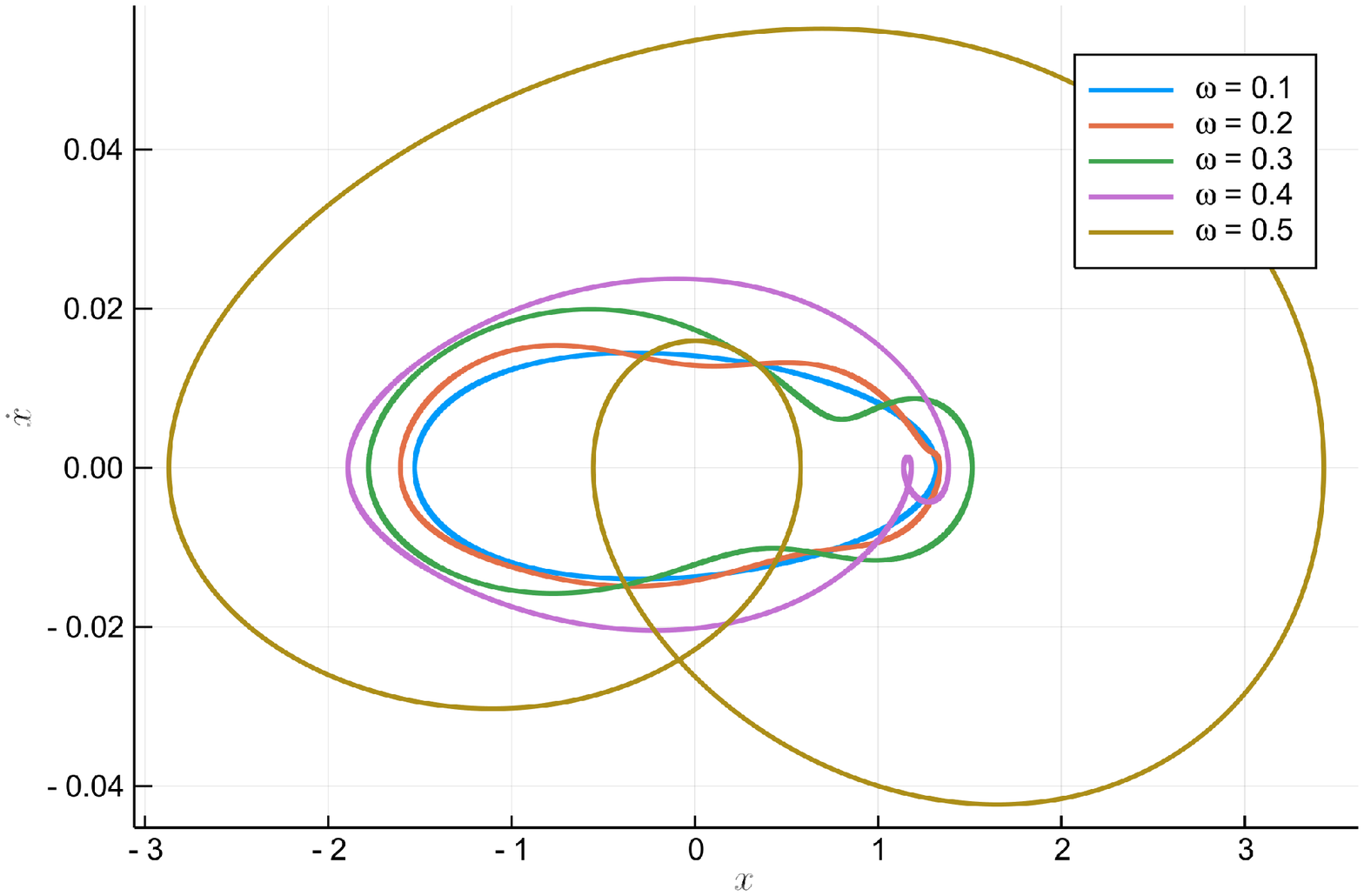}
	\end{subfigure}%
	~ 
	\begin{subfigure}[t]{0.5\textwidth}
		\centering
		\includegraphics[height=2in]{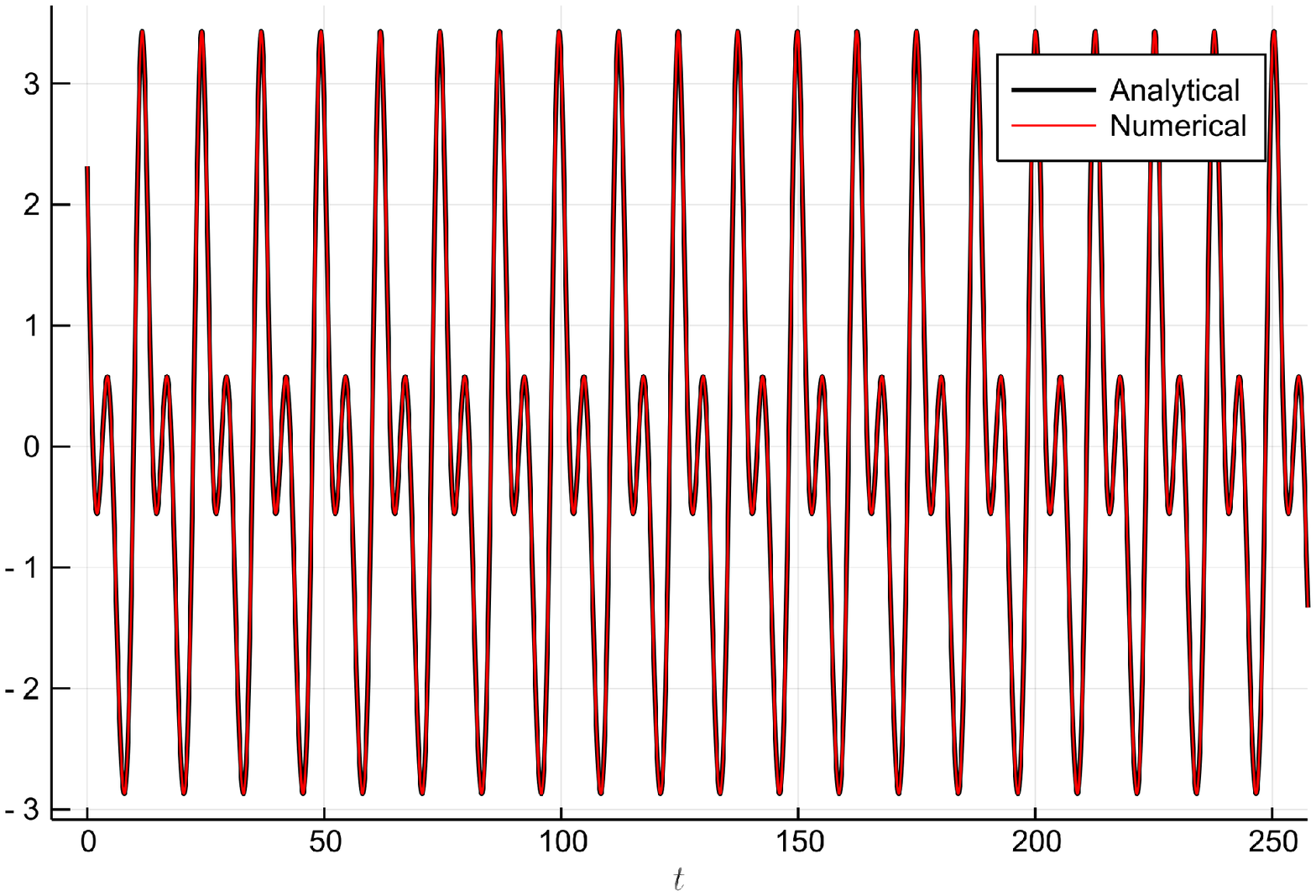}
	\end{subfigure}
	\caption{Quadratic model: (a)Phase space projections of the real parts of the displacement $x$ and the velocity $\dot{x}$ for the parameters $(k,h,f, g, \epsilon) = (1., 0.1, 0.8, 0, 0.1)$. (b) Time plots of the real parts of the analytical (black) and numerical (red) solutions with parameters as in (a) and $\omega = 0.5$.}\label{fig19}
\end{figure*}

Let us begin with the effect that the forcing frequency $\omega$ has on the periodic solutions of the quadratic nonlinear  model. Keeping our parameters fixed at $(k,h,f, g,\omega, \epsilon) = (1., 0.1, 0.8, 0., 0.1)$ we first plot in Fig. \ref{fig19} the solutions of the equation with the quadratic nonlinearity. 

Note in Fig. \ref{fig19}(a) that the projections of the periodic orbits on the $(x,\dot{x})$ plane, for various values of $\omega$, have a ``cyclical'' form and become more complex as we vary the forcing frequency. In Fig. \ref{fig19}(b) we fix $\omega = 0.5$ and plot the corresponding solution as a function of time, superimposing the analytical (black) and the numerical (red) results to demonstrate the accuracy of our computations. Clearly, the convoluted form of the nonlinear oscillations are in stark contrast with the corresponding phase space plots of the associated linear model which are simple ellipses describing harmonic motion. 

\begin{figure*}[t]
	\centering
	\begin{subfigure}[t]{0.5\textwidth}
		\centering
		\includegraphics[height=2in]{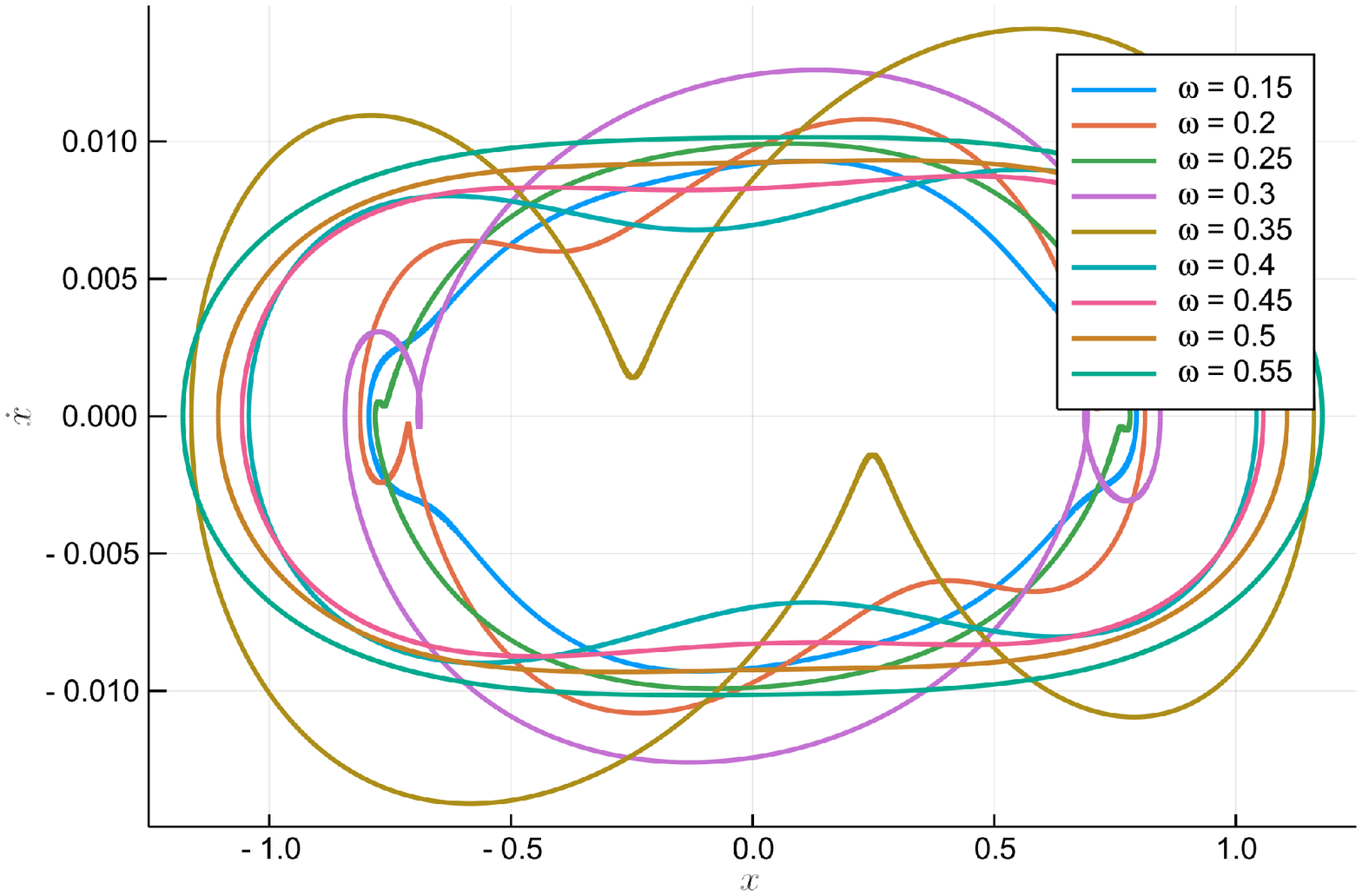}
	\end{subfigure}%
	~ 
	\begin{subfigure}[t]{0.5\textwidth}
		\centering
		\includegraphics[height=2in]{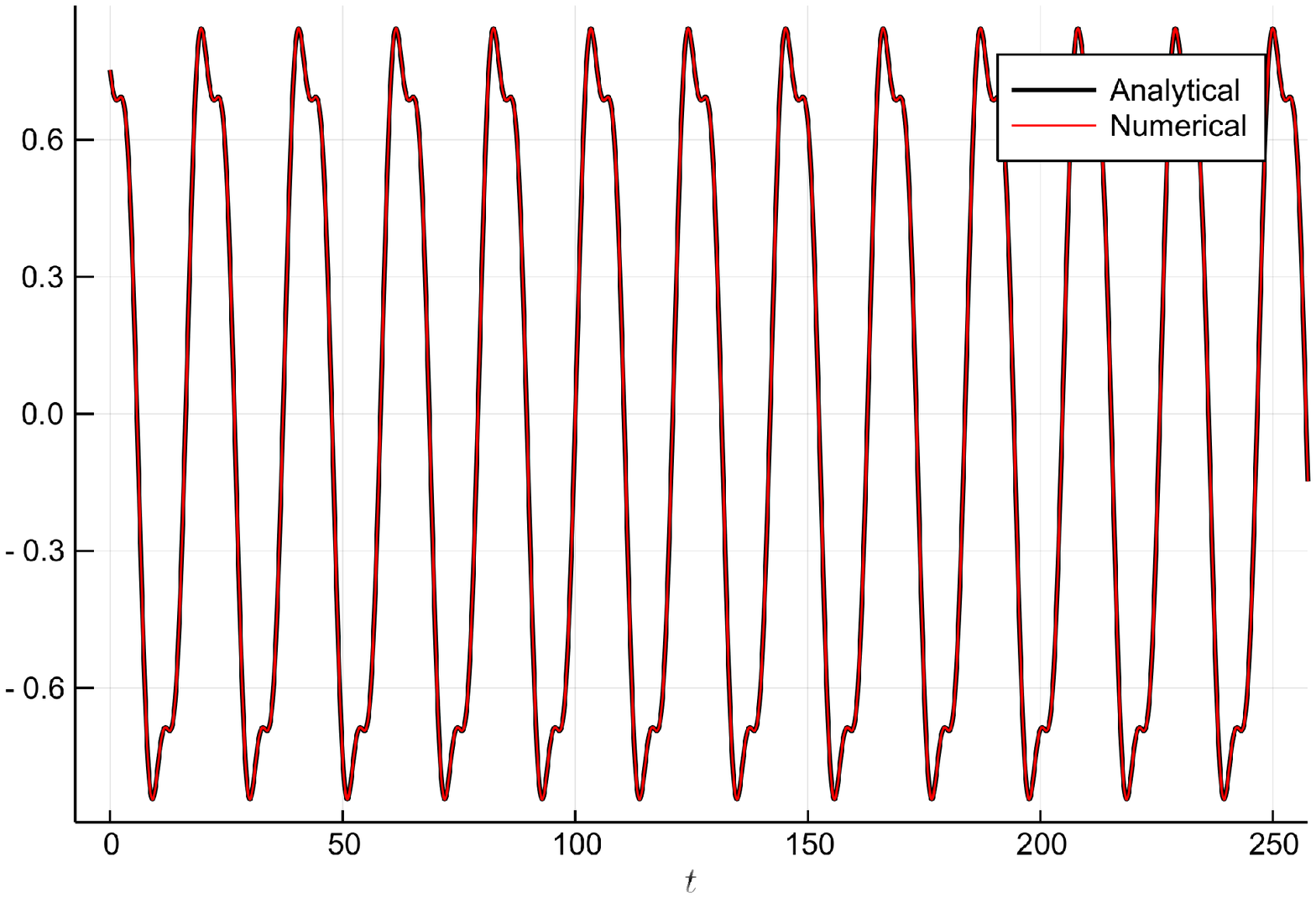}
	\end{subfigure}
	\caption{Cubic model: (a)Phase space projections of the real parts of the displacement $x$ and the velocity $\dot{x}$ for the parameters $(k,h,f, g, \epsilon) = (1., 0.1, 0.8, 0, 0.05)$. (b) Time plots of the real parts of the analytical (black) and numerical (red) solutions with parameters as in (a) and $\omega = 0.3$.}\label{fig20}
\end{figure*}

We now repeat the above study for the periodic solutions of the cubic model. Specifically, in Fig. \ref{fig20}(a) we set $(k,h,f, g, \epsilon) = (1., 0.1, 0.8, 0, 0.05)$ and vary the forcing frequency $\omega$, while in Fig. \ref{fig20}(b) we present the associated time plots of the real part of the displacement, for $\omega = 0.3$. Notice here also a complex morphology of the orbits associated with the cubic nonlinearity, similar to what we found in the case of the quadratic nonlinearity.

\begin{figure*}[t]
	\centering
	\begin{subfigure}[t]{0.5\textwidth}
		\centering
		\includegraphics[height=2in]{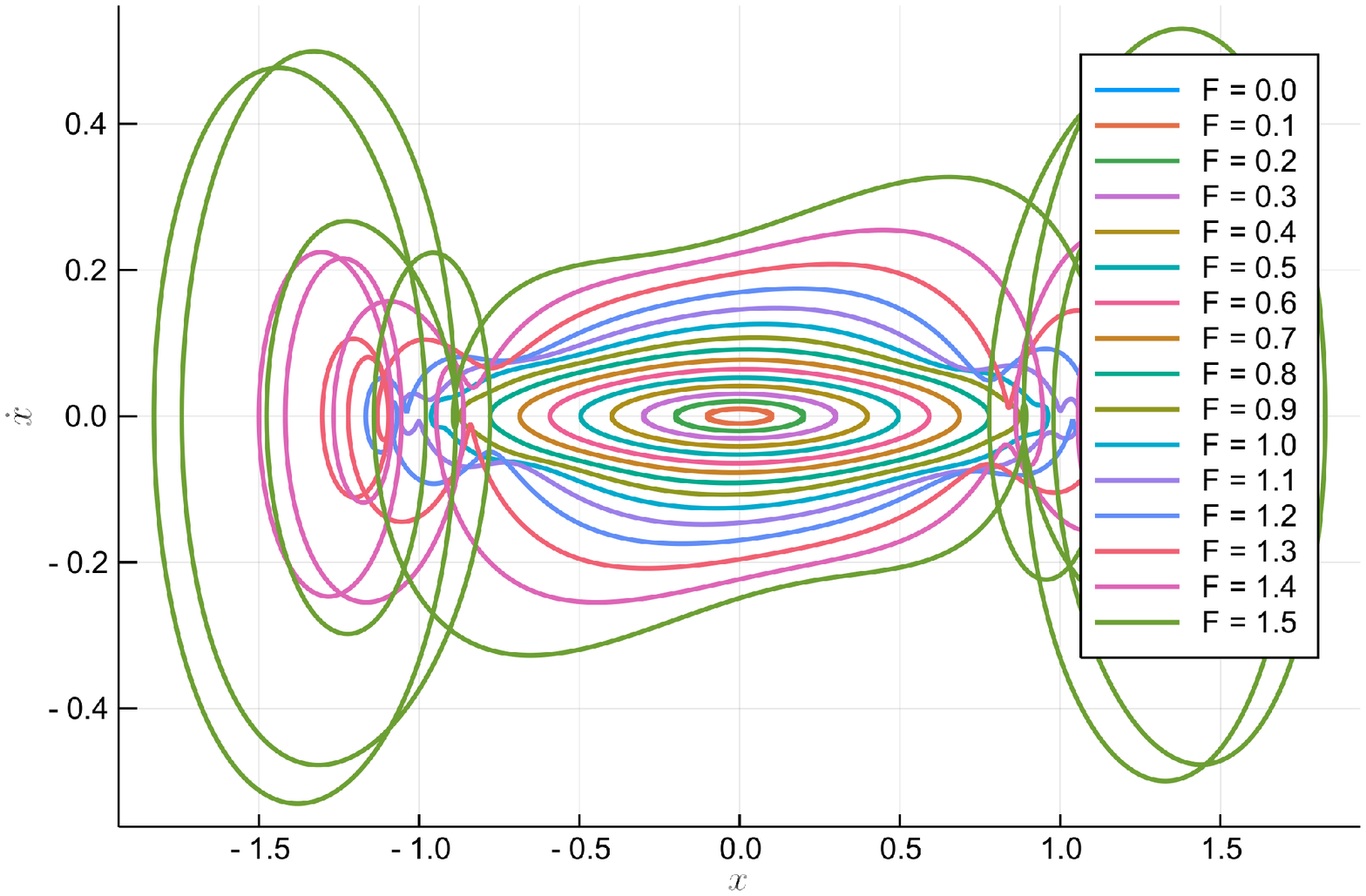}
	\end{subfigure}%
	~ 
	\begin{subfigure}[t]{0.5\textwidth}
		\centering
		\includegraphics[height=2in]{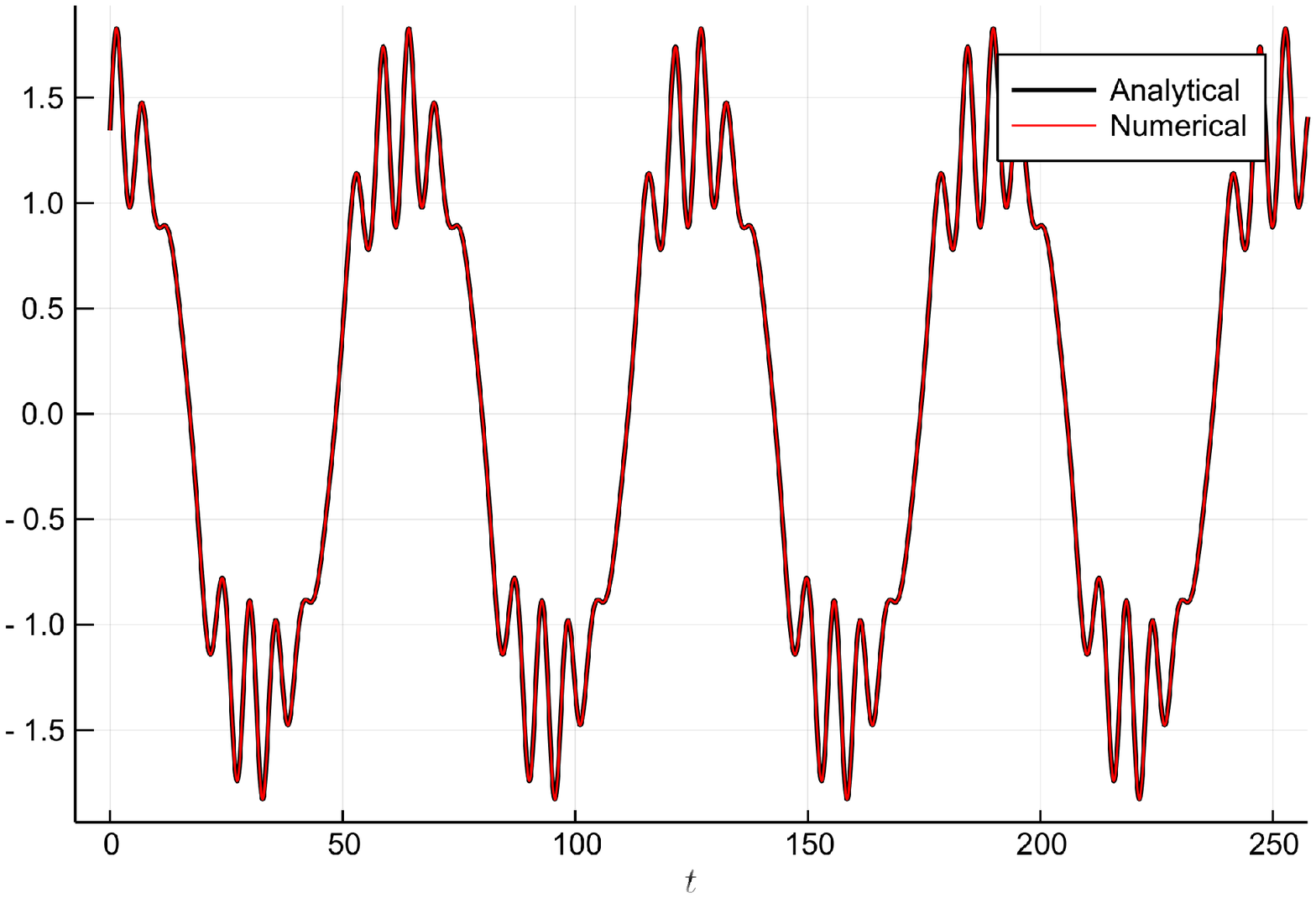}
	\end{subfigure}
	\caption{Forcing amplitude variation for the cubic model: (a)Phase space projections of the real parts of the displacement $x$ and the velocity $\dot{x}$ for the parameters $(k,h,f, g, \epsilon) = (1., 0.1, 0, 0.1, 0.05)$. (b) Time plots of the real parts of the analytical (black) and numerical (red) solutions with parameters as in (a) and $F = 0.3$.}\label{fig21}
\end{figure*}

Finally, let us investigate for the cubic nonlinearity model the effects of the driving amplitude on the periodic solutions, varying the value of $F$ (real) for $\omega = 0.1$. The results are even more remarkable and are presented here in Fig. \ref{fig21}. It is evident that the variation of the driving amplitude in this case results in an even greater increase of complexity in the morphology of the periodic solutions.

\subsubsection{Amplitude and phase response}

As we have seen earlier in the paper, a related important feature of the nonlinear models is the behavior of the amplitudes and phases of the periodic attractor as functions of the driving frequency, since more frequencies are present in the solution due to the nonlinearity of the equation of motion. In particular, we are interested in both \textit{amplitude and phase response curves} and their dependence on the parameters of the system. 

We recall that the amplitude response curve describes the ratio of the amplitude of the nonlinear attractor to the amplitude of the driving force, with respect to the frequency ratio $\omega/\omega_1$. Similarly, the frequency response curve describes the phase delay angle of the forced oscillation as a function of $\omega/\omega_1$.  As we have done for the quadratic nonlinearity model, we also compute here the amplitude and phase response curves of the present model using the corresponding quantities obtained from our analytical solutions.

\begin{figure}[H]
	\centering
	\includegraphics[width = 0.65\textwidth]{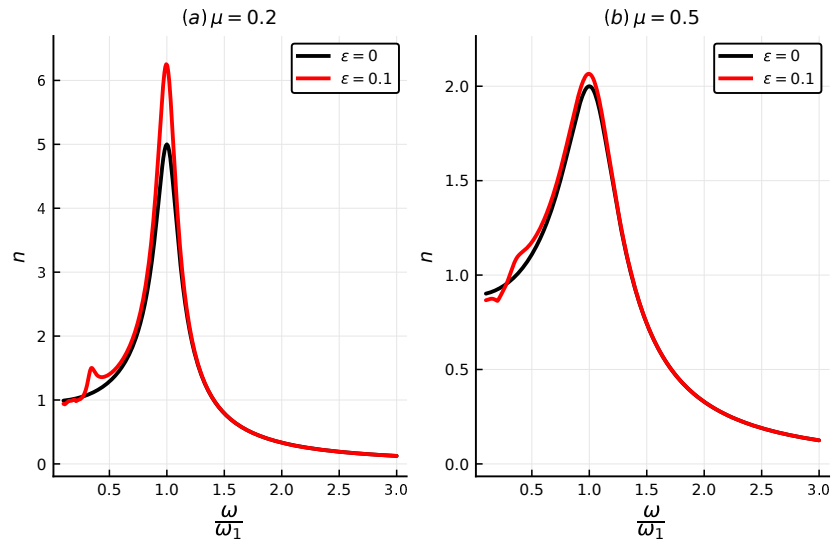}
	\caption{Magnification factor $n$ (amplitude response) with respect to the frequency ratio $\omega/\omega_1$ for the system with cubic nonlinearity.}	\label{fig22}
\end{figure}

\begin{figure}[H]
	\centering
	\includegraphics[width = 0.65\textwidth]{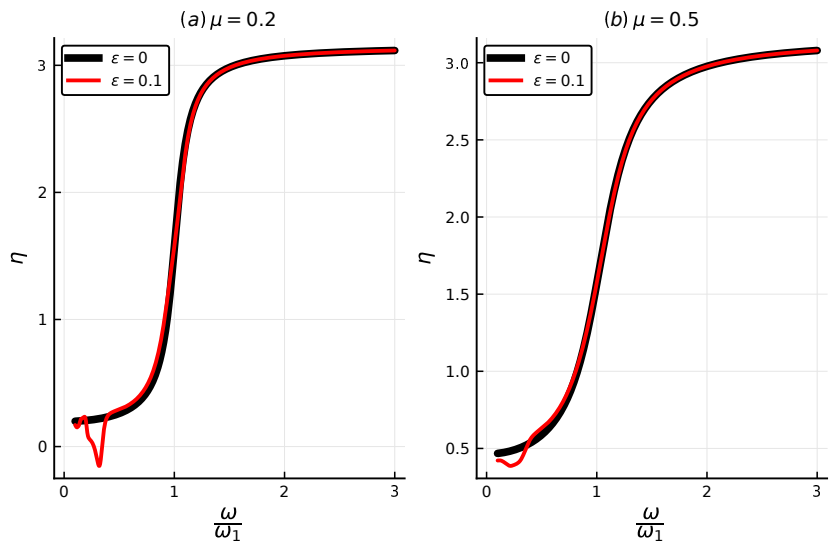}
	\caption{Phase response with respect to the frequency ratio $\omega/\omega_1$ for the system with cubic nonlinearity ($\epsilon=0.1$), superimposed with the associated response of the linear system ($\epsilon=0$).}	\label{fig23}
\end{figure}

In Fig. \ref{fig22} and Fig. \ref{fig23} respectively we present the amplitude response curve and the phase response curve of the cubic nonlinear model for $\mu=0.2, 0.5$ and $\epsilon = 0,0.1$, keeping $k=1$ fixed. Evidently, the results are qualitatively similar with what we found in the case of the quadratic nonlinearity.


\section{Reid's hysteretic oscillator with periodic forcing}

As we have seen in the previous sections there appear to be important difficulties when one tries to study numerically Bishop's complex hysteretic oscillator model, as introduced in \cite{bishop1955treatment}. As long as the model is kept linear, its general solution can be obtained analytically even in the case of $n$ such coupled oscillators \cite{bishop1956general}. However, when one tries to evaluate these solutions numerically, one finds that unavoidable errors arise, which cause these solutions to diverge to infinity after a time interval that depends on the parameters of the problem and the precision of the computations. These errors are also found to persist in very similar fashion, in the nonlinear versions of the model as well. Still, if one is content with studying the model for a specified time interval, periodic solutions are obtained which are in excellent agreement with the analytical ones, while the dynamics of the nonlinear model has features that are analogous to what one expects from nonlinear oscillators with viscous damping.

Hysteretic damping models, of course, have a long history, and there is a lively debate in the literature about their comparative virtues and shortcomings. Thus, we decided to include in this work a preliminary investigation of another popular hysteretic (in the sense of ``frequency--independent'') damping model, initially proposed by Reid \cite{reid1956free}. It is expressed in terms of a {\it real} (as opposed to complex) differential equation and hence its solutions are directly physically interpretable.

The properties of Reid's oscillator have been thoroughly analyzed in the literature  \cite{caughey1970free, caughey1976stability}. Its characteristic through-zero hysteretic loop shape in the force-displacement plot, having a bow-tie shape \cite{spitas2009continuous}, is clearly non-physical. Yet, for a wide range of engineering materials, the hysteresis loops are extremely narrow to the extent that their exact shape is neither noticeable nor significant for characterizing the dynamical response. Furthermore, various modifications of the model have also been proposed (e.g. in Refs. \cite{muravskii2004frequency, spitas2009continuous, liu2006reid}) to make it more realistic. Its main drawback appears to be its discontinuity at the points of stress--strain reversal, but it keeps drawing the attention of the scientific community mainly because of its simple form and its frequency independence property. In this section, we present a first study of its very interesting dynamical behavior in the presence of a cubic nonlinearity term in the stiffness potential.

Specifically, the evolution of the one degree of freedom Reid oscillator under periodic forcing is described by the following ``quasi--linear'' real differential equation:
\begin{equation}\label{reid1}
M\ddot{x} + c \abs{\dfrac{x}{\dot{x}}} \dot{x}  + k x = M\ddot{x} + k x \left(1 + \dfrac{c}{k} \text{sgn}(x \dot{x})\right) =  F \sin{\omega t},
\end{equation}
where $\text{sgn}(\cdot)$ is the sign function, $x$ denotes the particle displacement from equilibrium, $c$ is the damping coefficient, and $k$ quantifies the (linear) stiffness.

The oscillator described of Eqn. (\ref{reid1}) is frequency independent and yields work per cycle that is proportional the squared amplitude. Interestingly, it possesses stable periodic solutions (attractors), with frequency $\omega$, that are verified numerically for arbitrarily long times and are free from the spurious computational errors that plague Bishop's model. Here, we are interested in examining the dynamical effects of adding to the model a nonlinear cubic stiffness term, as we did for Bishop's oscillator, and compare the results. In this case, Reid's equation takes the form:

\begin{equation}\label{reid2}
M\ddot{x} + c \abs{\dfrac{x}{\dot{x}}} \dot{x}  + k x + \epsilon x^3 =  f \sin{\omega t}.
\end{equation}

\begin{figure}[H]
	\centering
	\includegraphics[scale=0.65]{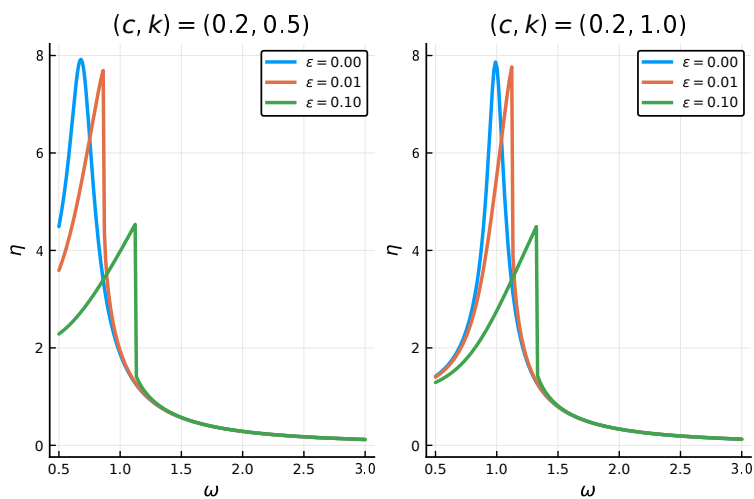}
	\caption{Amplitude response curves of Reid's model for varying $\epsilon\in \left\{0,0.01,0.1\right\}$. We set the damping coefficient to $c=0.2$ and use two different values for the linear stiffness, namely $k \in \left\{0.5, 1\right\}$.
		\label{fig_r1}}
\end{figure}

Although the original model of Eqn. (\ref{reid1}) is strictly speaking nonlinear, due to the sign--function term, this has limited implications for the dynamical behavior of the system. Note, for example, the differences in the system's amplitude response due to the cubic nonlinearity, as depicted in Figure \ref{fig_r1}. Although for high forcing frequencies the two systems behave almost indistinguishably, for low frequencies the model of Eqn. (\ref{reid2}) exhibits secondary resonances. Moreover, the peaks of the curves in Fig. \ref{fig_r1} are influenced by the nonlinearity, as higher values of $\epsilon$ lead to peaks at higher values of $\omega$, while we also observe that the higher the stiffness the lower the magnification factor. 

\begin{figure}[H]
	\centering
	\includegraphics[scale=0.65]{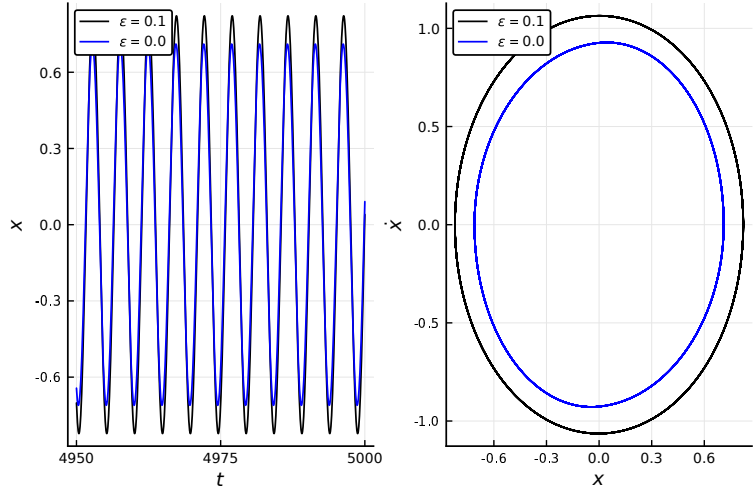}
	\caption{Periodic solutions for the two Reid's models, corresponding to the parameters $\left(c,k,f,\omega\right) = \left( 0.2, 0.3, 1., 1.3\right)$. Time plots (left) and phase plots (right), with $\epsilon \in \left\{0, 0.1\right\}$. \label{fig_r2}}
\end{figure}

For relatively high values of the damping (e.g. $c>0.1$), both systems are attracted to a stable periodic solution with period $T^{\star}=\frac{2\pi}{\omega}$, while, as we have already emphasized, both the original and extended models do not suffer from numerical instabilities of the type we encountered in the previous sections. In Fig. \ref{fig_r2} we present the periodic solutions for the two Reid's models, corresponding to the parameters $\left(c,k,f,\omega\right) = \left( 0.2, 0.3, 1., 1.3\right)$.

\begin{figure}[H]
	\centering
	\includegraphics[keepaspectratio, width=0.65\textwidth]{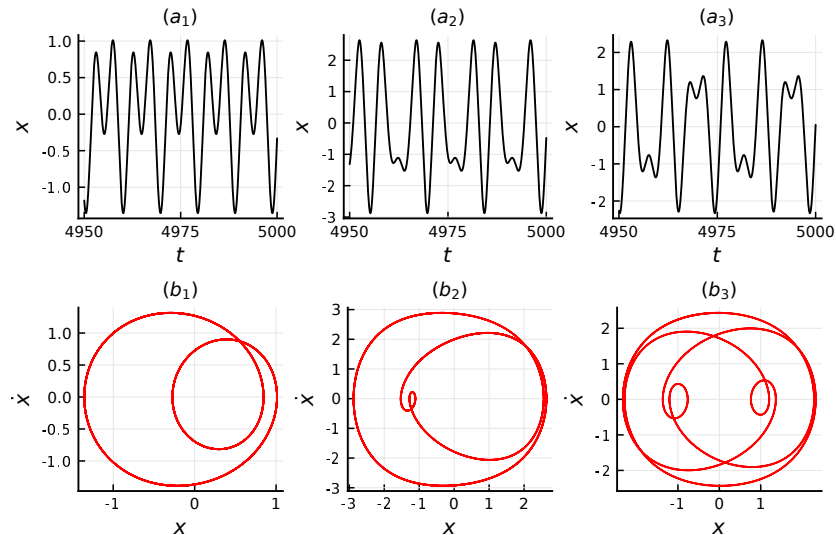}
	\caption{Time $(t,x(t))$ and phase $(x,\dot{x})$ plots of periodic solutions for the modified Reid's model, with parameters $\left(c,k,f,\omega,\epsilon \right) = \left(0.01, 0.3, 1.1, 1.3, 0.1\right)$. The associated periods are $2T^{\star}, 3T^{\star},5T^{\star}$ for the first, second and third column, respectively. \label{fig_r3}}
\end{figure}

On the other hand, for smaller damping coefficients $c$, the model of Eqn. (\ref{reid1}) eventually settles onto the stable periodic orbit with period $T^{\star}$, independently of the initial conditions $\left(x(0),\dot{x}(0)\right)$. This not the case, however, for the Reid's model of Eqn. (\ref{reid2}) with the cubic nonlinearity. As expected from analogous nonlinear models with viscous damping, for low values of $c$, one typically observes the emergence of stable periodic orbits with period different than $T^{\star}$. Such interesting stable solutions are presented here in Fig. \ref{fig_r3}, for parameters $\left(c,k,f,\omega,\epsilon \right) = \left(0.01, 0.3, 1.1, 1.3, 0.1\right)$. The associated periods are $2T^{\star}, 3T^{\star},5T^{\star}$ for the first, second and third column, respectively.

This is a remarkable observation that places Reid's nonlinear model in the wider framework of nonlinear oscillators known in the literature \cite{wiggins2003introduction}. In fact, when we plot the projections of Reid's modified model on the $x,dx/dt$ phase plane, we find a wealth of attracting periodic orbits, whose periods are multiples of $T^{\star}$ and can be reached, for the same parameter values, by merely changing the initial conditions! 

\begin{figure}[H]
	\centering
	\includegraphics[keepaspectratio, width=0.65\textwidth]{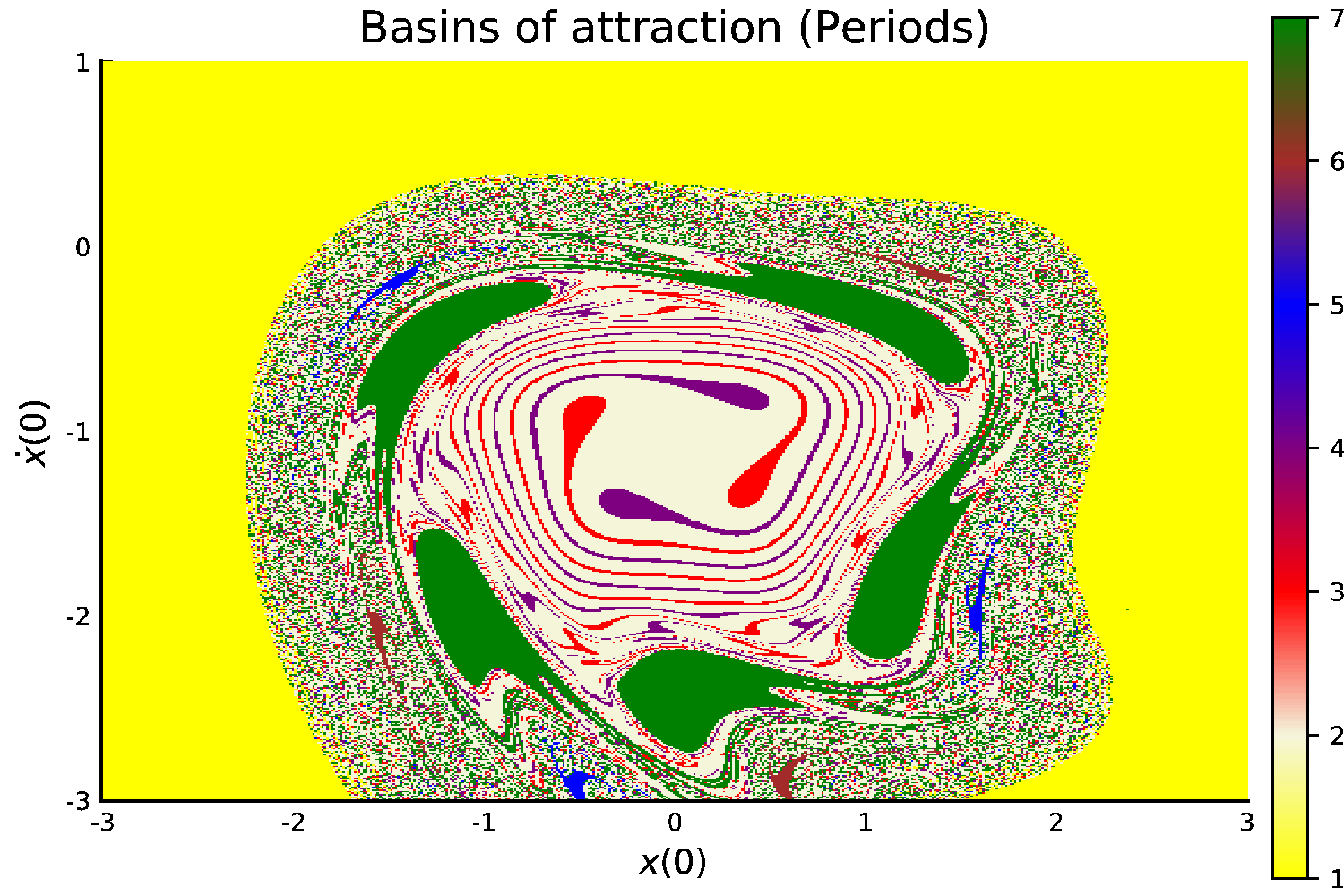}
	\caption{Basins of attraction of the modified Reid's model for parameters $\left(c,k,f,\omega,\epsilon \right) = \left(0.01, 0.3, 1.1, 1.3, 0.1\right)$. We have used the coding 1 and 2 for orbits with period $T = T^{\star}$, 3 and 4 for $T = 2 T^{\star}$,  5 and 6 for $T = 3T^{\star}$ and 7 for the $T =5 T^{\star}$. \label{fig_r4}}
\end{figure}

In fact, we investigate here pictorially this fascinating phenomenon of coexisting attractors by plotting their {\it basins of attraction} in the $x,dx/dt$ phase plane as follows: We select a grid of, say, $500 \times 500$ equally spaced initial conditions $(x(0),\dot{x}(0)) \in (-3,3) \times (-3,1)$, and color each one of them according to the attractor to which they lead as $t\rightarrow \infty$, for the same choice parameters. We present a preliminary study of these basins of attraction in Fig. \ref{fig_r4}, using different colors for the basins of attraction associated with the various periodic orbits. Specifically, we use yellow $(1)$ and beige $(2)$ for the orbits with periods $T^{\star}$, red $(3)$ and purple $(4)$ for the orbits with period $2T^{\star}$, blue $(5)$ and brown $(6)$ for the orbits with period $3T^{\star}$ and finally green $(7)$ for the orbit with period $5T^{\star}$ (the numbers in the previous parentheses refer to the colorbar coding of Fig. \ref{fig_r4}).

In conclusion, we have studied the dynamics of the hysteretic Reid's model in the presence of a cubic nonlinearity and have observed significant differences in its amplitude response with respect to the original oscillator. Depending on the magnitude of the nonlinearity parameter $\epsilon$, a remarkable coexistence of periodic attractors is revealed, whose periods are multiples of the driving period. This phenomenon of multistability is reminiscent of what one finds in nonlinear models with viscous damping, but has distinct features of its own. Observe its fractal--looking structure in the outer part of the plotted region in Fig. \ref{fig_r4}, where the system exhibits high sensitivity in the choice of initial conditions. The occurrence of such phenomena in coupled {\it systems} of Reid nonlinear oscillators needs to be further studied, as it might constitute an important first step in the analysis of materials with nonlinear behavior, such as multistable mechanical meta-materials \cite{rafsanjani2015snapping, yang2019multi, zhu2019bio}.


\section{Conclusions and discussion}

Oscillator models that combine nonlinearity with different forms of damping have lately attracted great interest, notably in cases where the damping involves time--dependent coefficients \cite{OlejAwre2018}, or is nonlinear due to stick--slip conditions \cite{OlejAwre2013}. In this paper the implications of the use of the Bishop and Reid hysteretic models with regard to the dynamical behavior and numerical accuracy and stability of single degree of freedom systems with quadratic and cubic stiffness non-linearity, corresponding to a composite meta-material, was studied under steady-state conditions. The key findings were:

1) We first solved the periodically driven Bishop's oscillator, in the presence of small stiffness nonlinearities, using perturbation theory in powers of a small parameter $\epsilon$ to obtain the periodic solutions in terms of a convergent Fourier series expansion. Our analytical findings allowed us to verify the validity of the numerical results and demonstrate that the amplitude of the periodic solutions increases significantly in magnitude and complexity with growing $\epsilon$.

2) We then extended our study and performed a comparative analysis using a different hysteretic model introduced by Reid \cite{reid1956free}, combining periodic forcing with weakly nonlinear stiffness terms. Reid's model is expressed in terms of a real differential equation that might be considered more ``physical'' than Bishop's complex model. Indeed, solving numerically Reid's model we always obtained stable periodic solutions as its attractors. More importantly, we discovered some remarkable multistability phenomena involving coexisting periodic attractors in the system's phase space, whose frequencies are integer multiples of the driving frequency.

The methods and results reported in this work may be used to aid in the selection and mathematical formulation of appropriate, realistic, computationally robust and numerically stable engineering models for nonlinear engineering materials. In this direction, we suggest that {\it arrays} of Reid nonlinear oscillators coupled by nearest neighbor interactions may have concrete applications. Preliminary results indicate that such arrays, when driven at one end by periodic forcing whose frequency lies outside the linear spectrum, can transmit low energy waves through the array in a highly controlled fashion. However, when the driving amplitude exceeds a certain threshold, a high energy wave is suddenly excited causing the normal operation of the system to break down \cite{BouKalSpi2020}.

\section*{Acknowledgments}
We are grateful to the reviewers for helping us improve our statement of the objectives, the outcomes and the conclusions, as well as updating our list of references. We acknowledge partial support for this work by funds from the Ministry of Education and Science of Kazakhstan, in the context of the project VSAT (2018-2020) and the Nazarbayev University internal grant HYST (2018-2021).

\bibliographystyle{asmems4}

\bibliography{asme2e}


\end{document}